\documentclass[a4paper,oneside,article,english,10pt]{memoir}

\usepackage{tex/preamble}

\newcommand\fag{Goodness in $n$-angulated categories}

\begin{document}

\newcommand\forfatter{Sebastian H. Martensen}
\newcommand\stdnr{(201506958)}
\newcommand\dato{\today}

{\centering
{\scshape\Large \fag}\par\vspace{0.25 cm}
{ \forfatter}\par\vspace{-0.2cm}
}
\bigskip


\newcommand\blfootnote[1]{%
  \begingroup
  \renewcommand\thefootnote{}\footnote{#1}%
  \addtocounter{footnote}{-1}%
  \endgroup
}

\makeatletter

\renewcommand\@makefntext[1]{%
\setlength\parindent{3pt}%
\noindent
\mbox{\@thefnmark~}{#1}}

\makeatother

\begin{abstract}
	\noindent We generalise the notions of good, middling good, and Verdier good morphisms of distinguished triangles in triangulated categories, first introduced by Neeman, to the setting of $n$-angulated categories, introduced in Geiss, Keller, and Oppermann. We then prove that all morphisms of $n$-angles in an $(n-2)$-cluster tilting $n$-angulated category are middling good for $n>3$.
	\blfootnote{%
		\emph{Date:} \today \par
		2020\emph{ Mathematics Subject Classification.} Primary 18G80.\par
	}
\end{abstract}

\chapter{Introduction}

Triangulated categories were introduced independently by Puppe and Verdier in 1962 and 1963 and, today, they turn up in the settings of representation theory, algebraic geometry, and algebraic topology to name a few. However, even in a model-free setting, the triangulation of a category is a rich structure. The Adams spectral sequence, which is classically defined over the stable homotopy category and calculates stable homotopy groups from mod $p$ cohomology operations, extends to general triangulated categories that are equipped with a \emph{projective} or \emph{injective class} (see \cite{christensen98}). This more general Adams spectral sequence is dependent on the chosen class and, given an appropriate class, it will for example specialise to the Adams--Novikov spectral sequence in the category of spectra and, in the derived category of a differential graded algebra, it can specialise to a spectral sequence that converges to sets of morphisms.

Christensen and Frankland found (see \cite{cf20f}) that, by choosing suitable fill-in morphisms between quotients in Adams towers, one might be able to generalize the Moss pairing to the general Adams spectral sequence in triangulated categories. 
%
%
Of interest in this context are fill-ins which are both \emph{good} and \emph{Verdier good} (see Definitions 3.1 and 3.4 below).

In short, a morphism of distinguished triangle is called \emph{good} if its mapping cone (see Definition \ref{mappingconedef} below) is a distinguished triangle, it is called \emph{Verdier good} if the third morphism between objects factors through two specific octahedra, and furthermore we call a morphism \emph{middling good} if it extends to a $(4\times4)$-diagram wherein the first three rows and columns are all distinguished triangles. Starting from a commutative square, it turns out that all three conditions are almost entirely dependent on the choice of fill-in map between the cone objects. In \cite{frankland20} and \cite{neeman91}, certain relationships between these concepts are found and certain classes of morphisms of distinguished triangles are found to be good if and only if they are Verdier good. It is certain, however, that both goodness and Verdier goodness are stronger conditions than middling goodness. Something to take note of, is that even with all the extra structure inherent to the derived category $D(\Z)$, which is an \emph{algebraic triangulated category}, there are simple examples of morphisms which are not middling good; take for instance the following diagram arising from a mapping of the embedding of the exact sequence
\begin{center}
	\begin{tikzpicture}	
		\diagram{d}{2em}{3em}{
		  0 & \Z & \Z & \Z/n & 0\\ 
		};

		\path[->, font = \scriptsize, auto]
		  (d-1-2) edge node{$n$} (d-1-3)
		  (d-1-3) edge node{$q_n$} (d-1-4)
		  (d-1-1) edge node{$$} (d-1-2)
		  (d-1-4) edge node{$$} (d-1-5)
		  ;
	\end{tikzpicture}
\end{center}
in degree 0 to its left rotation:

\begin{center}
	\begin{tikzpicture}	
		\diagram{d}{2em}{3em}{
		  \Z[0] & \Z[0] & \Z/n[0] & \Z[1]\\
		  \Z[0] & \Z/n[0] & \Z[1] & \Z[1].\\ 
		};

		\path[->, font = \scriptsize, auto]
		  (d-1-1) edge node{$n$} (d-1-2)
		  (d-1-2) edge node{$q_n$} (d-1-3)
		  (d-1-3) edge node{$\eps_n$} (d-1-4)
		  
		  (d-2-1) edge node{$q_n$} (d-2-2)
		  (d-2-2) edge node{$\eps_n$} (d-2-3)
		  (d-2-3) edge node{$-n$} (d-2-4)
		  
		  (d-1-1) edge node{$0$} (d-2-1)
		  (d-1-2) edge node{$0$} (d-2-2)
		  (d-1-3) edge node{$\eps_n$} (d-2-3)
		  (d-1-4) edge node{$0$} (d-2-4);
	\end{tikzpicture}
\end{center}
One can show explicitly that this morphism does not extend to a $(4\times 4)$-diagram or one can appeal to \cite[Proposition 7.13]{frankland20}.

In the setting of $n$-angulated categories, as introduced by Geiss, Keller, and Oppermann in \cite{gko}, even simple questions regarding the relationship of these three goodness criteria are difficult to answer; the extra commutativity relations and the distance between objects in an $n$-angle become obstructions to the tools and methods readily available in the triangulated setting. So far, the most concrete examples of $n$-angulated categories are those arising from $(n-2)$-cluster tilting subcategories of triangulated categories described in \cite{gko} and the so-called exotic $n$-angulated categories described in \cite{bjt}. 
\smallskip 

What we shall see in this paper, is that in $(n-2)$-cluster tilting $n$-angulated categories, not only do all commutative squares extend to an $(n+1)\times(n+1)$-diagram wherein the first $n$ rows and columns are $n$-angles, in fact, all morphisms of $n$-angles extend to such a diagram for $n>3$ (Theorem \ref{maintheorem}). Since $(n-2)$-cluster tilting $n$-angulated categories are algebraic $n$-angulated categories in the sense of \cite{jasso}, it is somewhat surprising that something similar does not hold for, say, $D(\Z)$. Naturally one might ask whether this holds true in general for $n$-angulated categories with $n>3$. This is, as of yet, unknown.\smallskip

We begin in Section 2 by recalling the standard definitions related to $n$-angulated categories. In Section 3.1 we recall the definitions and some of the results about the various goodness criteria, and in Section 3.2 we generalise them to the $n$-angulated setting. In Section 4.1 we recall the definition and $n$-angulation of $(n-2)$-cluster tilting categories and prove our main result, and in Section 4.2 we discuss the exotic case and open questions.\smallskip

\subsection*{Acknowledgements}
I would like thank Marius Thaule for the many discussions that resulted in this project and for all of the advice he has imparted along the way. I would also like to thank Petter Bergh and Martin Frankland for their helpful comments and interest in the project.


\chapter{Preliminary theory}

Let \cat C be an additive category equipped with an automorphism $\Sigma\colon\cat C\to\cat C$ and let $n$ an integer greater than or equal to three. 

\begin{definition}
	A diagram in \cat C of the form
\begin{center}
	\includegraphics{img/standalone/nAngle.tex}
\end{center}
	is called an \emph{$n$-$\Sigma$-sequence}. We may denote such a sequence by $A_\bullet$ and refer to the maps $\alpha_i$ as the \emph{structure maps} of the sequence. In particular, we may say that $\alpha_1$ is the \emph{base} of the $n$-$\Sigma$-sequence and that the $\alpha_i$ are the \emph{cones of} $\alpha_1$ for $i\in\{2,\dots,n\}$.
\end{definition}

\begin{definition}
	The \emph{left} and \emph{right rotations} of $A_\bullet$ are the $n$-$\Sigma$-sequences
	\begin{center}
		\includegraphics{img/standalone/leftRotation.tex}
	\end{center}
	and
	\begin{center}
		\includegraphics{img/standalone/rightRotation.tex}
	\end{center}
	respectively. We will borrow the notation of the shift functor from cochain complexes and thus denote by $A_\bullet[1]$ the left rotation of $A_\bullet$ and by $A_\bullet[-1]$ the right rotation. For an object $A$ of \cat C, a sequence of the form
	\begin{equation}\label{trivial}
	\begin{split}
		\includegraphics{img/standalone/trivialNAngle.tex}
	\end{split}
	\end{equation}
	and any of its rotations are called \emph{trivial} $n$-$\Sigma$-sequence. The $n$-$\Sigma$-sequence (\ref{trivial}), in particular, we denote by $(\Gamma A)_\bullet$.\smallskip

	An $n$-$\Sigma$-sequence $A_\bullet$ is \emph{exact} if the induced sequence 
	\begin{center}
		\includegraphics{img/standalone/homSigmaSeq.tex}
	\end{center}
	is exact as a cochain complex.
\end{definition}

\begin{definition}
	A \emph{morphism} of $n$-$\Sigma$-sequences $\varphi_\bullet\colon A_\bullet\to B_\bullet$ is an n-tuple of morphisms $\varphi_\bullet=(\varphi_1,\varphi_2,\dots,\varphi_n)$ constituting a commutative diagram
	\begin{center}
		\includegraphics{img/standalone/sigmaSeqMorphism.tex}
	\end{center}
The morphism $\varphi_\bullet$ is an \emph{isomorphism} if $\varphi_i$ is an isomorphism for each $i\in\{1,\dots,n\}$;  if for some $1\<i\<n$ the maps $\varphi_i$ and $\varphi_{i+1}$ are isomorphisms (identifying $\varphi_{n+1}\coloneqq\Sigma\varphi_1$) we say that $\varphi_\bullet$ is a \emph{weak isomorphism}.
\end{definition}

\begin{remark}
The collection of $n$-$\Sigma$-sequences over \cat C and the collection of morphisms between them constitutes a category. We may thus consider $\Gamma$ a functor from $\cat C$ to $n$-$\Sigma$-sequences taking an object to the trivial $n$-$\Sigma$-sequence (\ref{trivial}) and extending a morphism of objects to a morphism of $n$-$\Sigma$-sequences in the obvious way as illustrated below.
\begin{center}
	\includegraphics{img/standalone/actionOfGamma.tex}
\end{center}\vspace{-1.3cm}
\end{remark}

\begin{definition}\label{definitionnangulated}
Let $\mathscr N$ denote a collection of $n$-$\Sigma$-sequences in \cat C. The triple $(\cat C,\Sigma,\mathscr N)$ is an \emph{$n$-angulated category} if it satisfies the following four axioms:
\begin{enumerate}
	\item[(N1)]	
	\begin{enumerate}
		\item[(a)]	$\mathscr N$ is closed under direct sums, direct summands, and isomorphisms of $n$-$\Sigma$-sequences.
		\item[(b)]	For any object $A$ in \cat C, the trivial $n$-$\Sigma$-sequences $(\Gamma A)_\bullet$ is in $\mathscr N$.
		\item[(c)]	For each morphism $\alpha_1\colon A_1\to A_2$ in \cat C there exists an $n$-$\Sigma$-sequence in $\mathscr N$ with base $\alpha_1$.
	\end{enumerate}
	\item[(N2)]		An $n$-$\Sigma$-sequence belongs to $\mathscr N$ if and only if its left rotation belongs to $\mathscr N.$
	\item[(N3)]		Given that the solid part of the diagram\vspace{-.5em}
	\begin{center}
		\begin{tikzpicture}	
			\diagram{d}{2em}{3em}{
			  A_1 & A_2 & A_3 & \cdots & A_n & \Sigma A_1\\
			  B_1 & B_2 & B_3 & \cdots & B_n & \Sigma B_1,\\ 
			};
	
			\path[->, font = \scriptsize, auto]
			  (d-1-1) edge node{$\alpha_1$} (d-1-2)
			  (d-1-2) edge node{$\alpha_2$} (d-1-3)
			  (d-1-3) edge node{$\alpha_3$} (d-1-4)
			  (d-1-4) edge node{$\alpha_{n-1}$} (d-1-5)
			  (d-1-5) edge node{$\alpha_n$} (d-1-6)
			  
			  (d-2-1) edge node{$\beta_1$} (d-2-2)
			  (d-2-2) edge node{$\beta_2$} (d-2-3)
			  (d-2-3) edge node{$\beta_3$} (d-2-4)
			  (d-2-4) edge node{$\beta_{n-1}$} (d-2-5)
			  (d-2-5) edge node{$\beta_n$} (d-2-6)

			  (d-1-1) edge node{$\varphi_1$} (d-2-1)
			  (d-1-2) edge node{$\varphi_2$} (d-2-2)
			  (d-1-3) edge[densely dashed] node{$\varphi_3$} (d-2-3)
			  (d-1-5) edge[densely dashed] node{$\varphi_n$} (d-2-5)
			  (d-1-6) edge node{$\Sigma\varphi_1$} (d-2-6)
			  ;
\end{tikzpicture}
	\end{center}\vspace{-1em}
	 commutes and has rows in $\mathscr N$, the dotted morphisms exist such that $(\varphi_1,\varphi_2,\varphi_3,\dots,\varphi_n)$ is a morphism of $n$-$\Sigma$-sequence.
	\item[(N4)]		Given that the solid part of Figure \ref{octahedron} commutes
	\begin{figure}[h]
		\centering
		\includegraphics{img/standalone/octahedron.tex}
		\caption{An octahedron in an $n$-angulated category.}
		\label{octahedron}
	\end{figure}
	 and has rows in $\mathscr N$, the dotted morphisms exist such that each square commutes and such that the $n$-$\Sigma$-sequence of Figure \ref{assoc} belongs to $\mathscr N$.\vspace{-1.5em}
	
	\begin{center}
	\begin{figure}[h]
		\begin{tikzpicture}
			\node (1) at (-2,1.14){$A_{3}$};
			\node (2) at (0,1.14){$A_4\oplus B_3$};
			\node (3) at (3.5,1.14){$A_5\oplus B_4\oplus C_3$};
			\node (4) at (6.5,1.14){$A_6\oplus B_5\oplus C_4$};
			\node (5) at (-3,0){$\cdots$};
			\node (6) at (0,0){$A_n\oplus B_{n-1}\oplus C_{n-2}$};
			\node (7) at (3,0){$B_n\oplus C_{n-1}$};
			\node (8) at (5.7,0){$C_n$};
			\node (9) at (7.7,0){$\Sigma A_3$};

			\draw[->] (4)+(1.2,0) arc (90:-90:.25cm) --node[fill=white,font=\scriptsize]{$\mu_{4}$} +(-11.2,0) arc (90:270:.3cm) ;
	
			\path[->, font = \scriptsize, auto]
			  (1) edge node{
			  $\big(\begin{smallmatrix}
				\alpha_3\\
				\varphi_3
			\end{smallmatrix}\big)$
			} (2)
			  (2) edge node{
			  $\bigg(\begin{smallmatrix}
				-\alpha_4&0\\
				\varphi_4&-\beta_3\\
				\lambda_4&\psi_3
			  \end{smallmatrix}\bigg)$
			} (3)
			  (3) edge node{$\mu_3$} (4)
			  (5) edge node{$\mu_{n-3}$} (6)
			  (6) edge node{$\eta$} (7)
			  (7) edge node{
			  $(\begin{smallmatrix}
				\psi_{n}&\gamma_{n-1}
			  \end{smallmatrix})$} (8)
			  (8) edge node{$(\Sigma\alpha_2)\gamma_n$} (9)
			  ;
			  \node (a) at (0,-1) {$\mt{with}\quad\mu_i=\bigg(\begin{smallmatrix}
				-\alpha_{i+2}&0&0\\
				(-1)^i\varphi_{i+2}&-\beta_{i+1}&0\\
				\lambda_{i+2}&\psi_{i+1}&\gamma_{i}
			  \end{smallmatrix}\bigg)\quad \mt{and}$};
			  \node (b) at (5.2,-1) {$\eta=\big(\begin{smallmatrix}
				(-1)^{n}\varphi_{n}&-\beta_{n-1}&0\\
				\lambda_{n}&\psi_{n-1}&\gamma_{n-2}
			  \end{smallmatrix}\big).$};
\end{tikzpicture}
	\caption{The $n$-angle associated to the octahedron.}
	\label{assoc}
	\end{figure}
	\end{center}
\end{enumerate}\vspace{-.5cm}
	The collection $\mathscr N$ is called an \emph{$n$-angulation} of the category \cat C with respect to the autoequivalence $\Sigma$, and the elements of $\mathscr N$ are called \emph{$n$-angles}. A category 
$\cat T$ is triangulated if and only if it is 3-angulated, so one may refer to the distinguished triangles of \cat T as 3-angles should one choose.
\end{definition}

We will often refer to the diagram of Figure \ref{octahedron} as an \emph{octahedron} and to the $n$-angle of Figure \ref{assoc} as the \emph{associated $n$-angle}. We will also refer to the axioms by their name (N1)--(N4), and in the case of triangulated categories we may at times refer to the relevant axioms as (TR1)--(TR4).

\begin{remark}
Similarly to the triangulated setting, (N1)(ii), (N2), and (N3) ensures that for an $n$-angle $A_\bullet$ with structure morphisms $\alpha_\bullet$, the map $\alpha_i$ is both a weak cokernel of $\alpha_{i-1}$ and a weak kernel of $\alpha_{i+1}$. 
\end{remark}

The category of $n$-$\Sigma$-sequences, in, for example, an $n$-angulated category, shares some similarities with categories of cochain complexes over abelian categories and, as such, we can borrow some concepts from this theory. As an example, one could bring the concept of chain homotopies to the setting of $n$-$\Sigma$-sequences as in \cite{bjt} and thus we can think of an $n$-angle as being \emph{nullhomotopic} (or \emph{contractible}) in the sense of regular chain complexes. 

Of course, when we have a notion of homotopy, the next thing we want is a notion of mapping cone, ideally one that is homotopy invariant. This already exists in categories of cochain complexes, and the idea directly translates to the setting of $n$-$\Sigma$-sequences. 

\begin{definition}\label{mappingconedef}
	Let $\varphi_\bullet\colon A_\bullet\to B_\bullet$ be a morphism of $n$-$\Sigma$-sequences
	\begin{center}
		\begin{tikzpicture}	
			\diagram{d}{2em}{3em}{
			  A_1 & A_2 & A_3 & \cdots & A_n & \Sigma A_1\\
			  B_1 & B_2 & B_3 & \cdots & B_n & \Sigma B_1.\\
			};
	
			\path[->, font = \scriptsize, auto]
			  (d-1-1) edge node{$\alpha_1$} (d-1-2)
			  (d-1-2) edge node{$\alpha_2$} (d-1-3)
			  (d-1-3) edge node{$\alpha_3$} (d-1-4)
			  (d-1-4) edge node{$\alpha_{n-1}$} (d-1-5)
			  (d-1-5) edge node{$\alpha_n$} (d-1-6)
			  
			  (d-2-1) edge node{$\beta_1$} (d-2-2)
			  (d-2-2) edge node{$\beta_2$} (d-2-3)
			  (d-2-3) edge node{$\beta_3$} (d-2-4)
			  (d-2-4) edge node{$\beta_{n-1}$} (d-2-5)
			  (d-2-5) edge node{$\beta_n$} (d-2-6)
			  
			  (d-1-1) edge node{$\varphi_1$} (d-2-1)
			  (d-1-2) edge node{$\varphi_2$} (d-2-2)
			  (d-1-3) edge node{$\varphi_3$} (d-2-3)
			  (d-1-5) edge node{$\varphi_n$} (d-2-5)
			  (d-1-6) edge node{$\Sigma\varphi_1$} (d-2-6)
			  ;
\end{tikzpicture}
	\end{center}
	Then the mapping cone of $\varphi_\bullet$ is the $n$-$\Sigma$-sequence
	\begin{center}
		\begin{tikzpicture}	
			\diagram{d}{1.5em}{3.5em}{
			  A_2\oplus B_1 & A_3\oplus B_2 & A_4\oplus B_3  \\\\
			  \cdots & \Sigma A_1\oplus B_{n} & \Sigma A_2\oplus\Sigma B_1.\\ 
			};
	
			\path[->, font = \scriptsize, auto]
			  (d-1-1) edge node{$\begin{psmallmatrix}
			  	-\alpha_2 & 0\\
			  	\varphi_2 & \beta_1
			  \end{psmallmatrix}$} (d-1-2)
			  (d-1-2) edge node{$\begin{psmallmatrix}
			  	-\alpha_3 & 0\\
			  	\varphi_3 & \beta_2
			  \end{psmallmatrix}$} (d-1-3)
			  (d-3-1) edge node{$\begin{psmallmatrix}
			  	-\alpha_n & 0\\
			  	\varphi_n & \beta_{n-1}
			  \end{psmallmatrix}$} (d-3-2)
			  (d-3-2) edge node{$\begin{psmallmatrix}
			  	-\Sigma\alpha_1 & 0\\
			  	\Sigma\varphi_1 & \beta_n
			  \end{psmallmatrix}$} (d-3-3)
			  ;
			  \draw[->] (d-1-3)+(.8,0) arc (90:-90:.3) --node[fill=white,pos=.55,font=\scriptsize]{
			  $\begin{psmallmatrix}
			  	-\alpha_4 & 0\\
			 	\varphi_4 & \beta_3
			  \end{psmallmatrix}$} +(-7.2,0) arc (90:270:.48);
\end{tikzpicture}
	\end{center}
\end{definition}

However, a mapping cone should be the same type of object as the source and target of the morphism from which it was constructed and, while the mapping cone is an $n$-$\Sigma$-sequence by construction, if we start from two $n$-angles in an $n$-angulated category, the resulting mapping cone may not be an $n$-angle. This inspires the following definition.

\begin{definition}\label{ngoodmorphism}
	In a $n$-angulated category, a morphism of $n$-angles $$\varphi_\bullet\colon A_\bullet\to B_\bullet$$ is a \emph{good morphism} if the mapping cone is an $n$-angle.
\end{definition}

\begin{remark}\label{N4equivalence}
	The \dq{octahedral axiom} (N4) has a few variants. The original definition of $n$-angulated categories stated in \cite{gko} says the following:\smallskip

	For any two $n$-angles $A_\bullet$ and $B_\bullet$ and a solid commutative square such as the one below, there exist fill-in morphisms $\varphi_3,\dots,\varphi_n$ such that the full diagram commutes and such that $(\varphi_1,\varphi_2,\varphi_3,\dots,\varphi_n)$ is a good morphism of $n$-angles.
	\begin{center}
		\begin{tikzpicture}	
			\diagram{d}{2em}{3em}{
			  A_1 & A_2 & A_3 & \cdots & A_n & \Sigma A_1\\
			  B_1 & B_2 & B_3 & \cdots & B_n & \Sigma B_1.\\
			};
	
			\path[->, font = \scriptsize, auto]
			  (d-1-1) edge node{$\alpha_1$} (d-1-2)
			  (d-1-2) edge node{$\alpha_2$} (d-1-3)
			  (d-1-3) edge node{$\alpha_3$} (d-1-4)
			  (d-1-4) edge node{$\alpha_{n-1}$} (d-1-5)
			  (d-1-5) edge node{$\alpha_n$} (d-1-6)
			  
			  (d-2-1) edge node{$\beta_1$} (d-2-2)
			  (d-2-2) edge node{$\beta_2$} (d-2-3)
			  (d-2-3) edge node{$\beta_3$} (d-2-4)
			  (d-2-4) edge node{$\beta_{n-1}$} (d-2-5)
			  (d-2-5) edge node{$\beta_n$} (d-2-6)
			  
			  (d-1-1) edge node{$\varphi_1$} (d-2-1)
			  (d-1-2) edge node{$\varphi_2$} (d-2-2)
			  (d-1-3) edge[densely dashed] node{$\varphi_3$} (d-2-3)
			  (d-1-5) edge[densely dashed] node{$\varphi_n$} (d-2-5)
			  (d-1-6) edge node{$\Sigma\varphi_1$} (d-2-6)
			  ;
\end{tikzpicture}
	\end{center}
	It was proven in \cite{berghthaule} that, fixing the axioms (N1)--(N3), this axiom is equivalent to an axiom defined in said article, called (N4*), and it was later proven in \cite{arentz} that (N4*) is equivalent to the axiom (N4) stated in Definition \ref{definitionnangulated}.
\end{remark}

\chapter{Goodness criteria}

In this section we recall the various goodness criteria in the original triangulated setting and generalise their respective definitions to the $n$-angulated setting.

\section{Triangulated categories}

Definition \ref{ngoodmorphism} is a generalisation of the following definition, originally introduced by Neeman in \cite{neeman91}.
\begin{definition}\label{goodmorphism}
	In a triangulated category, a morphism of distinguished triangles
	\begin{equation}\label{morphism}
\begin{tikzpicture}[baseline=(current bounding box.center)]
			\diagram{d}{2em}{3em}{
			  A_1 & A_2 & A_3 & \Sigma A_1\\
			  B_1 & B_2 & B_3 & \Sigma B_1\\ 
			};
	
			\path[->, font = \scriptsize, auto]
			  (d-1-1) edge node{$\alpha_1$} (d-1-2)
			  (d-1-2) edge node{$\alpha_2$} (d-1-3)
			  (d-1-3) edge node{$\alpha_3$} (d-1-4)
			  
			  (d-2-1) edge node{$\beta_1$} (d-2-2)
			  (d-2-2) edge node{$\beta_2$} (d-2-3)
			  (d-2-3) edge node{$\beta_3$} (d-2-4)
			  
			  (d-1-1) edge node{$\varphi_1$} (d-2-1)
			  (d-1-2) edge node{$\varphi_2$} (d-2-2)
			  (d-1-3) edge node{$\varphi_3$} (d-2-3)
			  (d-1-4) edge node{$\Sigma\varphi_1$} (d-2-4)
			  ;
\end{tikzpicture}
\end{equation}
	is a \emph{good morphism} if the mapping cone
	\begin{center}
		\begin{tikzpicture}
			\diagram{d}{2em}{4em}{
			  A_2\oplus B_1 & A_3 \oplus B_2 & \Sigma A_1\oplus B_3 &[1em] \Sigma A_2\oplus \Sigma B_1\\ 
			};
	
			\path[->, font = \scriptsize, auto]
			  (d-1-1) edge node{$
			  \Big(\begin{smallmatrix}
			  	-\alpha_2&0\\
			  	\varphi_{2}&\beta_1
			  \end{smallmatrix}\Big)
			  $} (d-1-2)
			  (d-1-2) edge node{$
			  \Big(\begin{smallmatrix}
			  	-\alpha_3&0\\
			  	\varphi_3&\beta_2
			  \end{smallmatrix}\Big)
			  $} (d-1-3)
			  (d-1-3) edge node{$
			  \Big(\begin{smallmatrix}
			  	-\Sigma\alpha_1&0\\
			  	\Sigma\varphi_1&\beta_3
			  \end{smallmatrix}\Big)
			  $} (d-1-4)
			  ;
\end{tikzpicture}
	\end{center}
	is a distinguished triangle.
\end{definition}

As was discussed in the case of $n$-angulated categories in Remark \ref{N4equivalence}, the octahedral axiom may be rephrased in terms of good morphisms. This is originally due to Neeman, the proof split between \cite[Theorem 1.8]{neeman91} and \cite[Proposition 1.4.6]{neeman01}.

\begin{theorem}[A. Neeman]
	Suppose \cat C is a pretriangulated category. The category \cat C is triangulated if and only if for each diagram as the one below with both rows consisting of distinguished triangles, a map $\varphi_3\colon A_3\to B_3$ exists making the diagram a good morphism. We may refer to $\varphi_3$ as a \emph{good fill-in}.
	\begin{equation}\label{diag:fill-in}
	\begin{aligned}
		\begin{tikzpicture}
			\diagram{d}{2em}{3em}{
			  A_1 & A_2 & A_3 & \Sigma A_1\\
			  B_1 & B_2 & B_3 & \Sigma B_1\\ 
			};
	
			\path[->, font = \scriptsize, auto]
			  (d-1-1) edge node{$\alpha_1$} (d-1-2)
			  (d-1-2) edge node{$\alpha_2$} (d-1-3)
			  (d-1-3) edge node{$\alpha_3$} (d-1-4)

			  (d-2-1) edge node{$\beta_1$} (d-2-2)
			  (d-2-2) edge node{$\beta_2$} (d-2-3)
			  (d-2-3) edge node{$\beta_3$} (d-2-4)

			  (d-1-1) edge node{$\varphi_1$} (d-2-1)
			  (d-1-2) edge node{$\varphi_2$} (d-2-2)
			  (d-1-3) edge[densely dashed] node{$\varphi_3$} (d-2-3)
			  (d-1-4) edge node{$\Sigma\varphi_1$} (d-2-4)
			  ;
		\end{tikzpicture}
		\end{aligned}
	\end{equation}
\end{theorem}

\begin{definition}
	A map of distinguished triangles (\ref{morphism}) is \emph{middling good} if it extends to a diagram
	\begin{center}
		\begin{tikzpicture}	
			\diagram{d}{2em}{3em}{
			  A_1 & A_2 & A_3 & \Sigma A_1\\
			  B_1 & B_2 & B_3 & \Sigma B_1\\
			  C_1 & C_2 & C_3 & \Sigma C_1\\
			  \Sigma A_1 & \Sigma A_2 & \Sigma A_3 & \Sigma^2 A_1 \\
			};
	
			\path[->, font = \scriptsize, auto]
			  (d-1-1) edge node{$\alpha_1$} (d-1-2)
			  (d-1-2) edge node{$\alpha_2$} (d-1-3)
			  (d-1-3) edge node{$\alpha_3$} (d-1-4)
			  
			  (d-2-1) edge node{$\beta_1$} (d-2-2)
			  (d-2-2) edge node{$\beta_2$} (d-2-3)
			  (d-2-3) edge node{$\beta_3$} (d-2-4)
			  
			  (d-3-1) edge[densely dashed] node{$$} (d-3-2)
			  (d-3-2) edge[densely dashed] node{$$} (d-3-3)
			  (d-3-3) edge[densely dashed] node{$$} (d-3-4)
			  
			  (d-4-1) edge node{$\Sigma\alpha_1$} (d-4-2)
			  (d-4-2) edge node{$\Sigma\alpha_2$} (d-4-3)
			  (d-4-3) edge node{$\Sigma\alpha_3$} (d-4-4)
			  
			  (d-1-1) edge node{$\varphi_1$} (d-2-1)
			  (d-1-2) edge node{$\varphi_2$} (d-2-2)
			  (d-1-3) edge node{$\varphi_3$} (d-2-3)
			  (d-1-4) edge node{$\Sigma\varphi_1$} (d-2-4)
			  
			  (d-2-1) edge node{$$} (d-3-1)
			  (d-2-2) edge node{$$} (d-3-2)
			  (d-2-3) edge node{$$} (d-3-3)
			  (d-2-4) edge node{$$} (d-3-4)
			  
			  (d-3-1) edge node{$$} (d-4-1)
			  (d-3-2) edge node{$$} (d-4-2)
			  (d-3-3) edge node{$$} (d-4-3)
			  (d-3-4) edge node{$$} (d-4-4)
			  ;

			  \path (d-3-3) -- (d-4-4) node[midway] (b) {$\Omega$};
\end{tikzpicture}
	\end{center}
	such that all squares commute, except for the square $\Omega$ which anticommutes, and such that the first three rows and columns are all distinguished triangles. We may refer to $\varphi_3$ as a \emph{middling good fill-in}.
\end{definition}

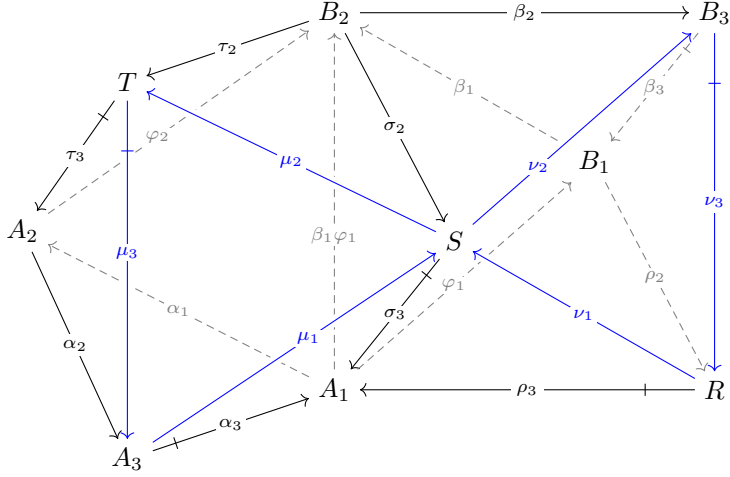
\begin{figure}[ht]
  \centering
  \begin{tikzpicture}[z=-10.5, scale = 2]
    \node (S) at (0,.2,1.25){$S$};
    \node (A1) at (-1.25,-1.25,0){$A_1$};
    \node (B2) at (-1.25,1.25,0){$B_2$};
    \node (T') at (1.25,-1.25,0){$R$};
    \node (B3) at (1.25,1.25,0){$B_3$};
    \node (B1) at (0,-.2,-1.25){$B_1$};
    \node (T) at (-2.15,1.25,1.25){$T$};
    \node (A3) at (-2.15,-1.25,1.25){$A_3$};
    \node (A2) at (-3.3,-.2,0){$A_2$};

    \path[->, font = \scriptsize, midway]
      (S) edge[susp] node[inner sep=1pt,fill=white]{$\sigma_3$} (A1)
      (A1) edge[densely dashed,gray] node[pos=.4,inner sep=1pt,fill=white]{$\beta_1\varphi_1$} (B2)
      (B2) edge node[inner sep=2pt,fill=white]{$\sigma_2$} (S)
      (B3) edge[susp,blue] node[inner sep=2pt,fill=white]{$\nu_3$} (T')
      (T') edge[susp] node[inner sep=2pt,fill=white]{$\rho_3$} (A1)
      (B1) edge[densely dashed,gray] node[inner sep=1pt,fill=white]{$\rho_2$} (T')
      (A1) edge[densely dashed,gray] node[pos=.45,inner sep=1pt,fill=white]{$\varphi_1$} (B1)
      (B2) edge node[inner sep=2pt,fill=white]{$\beta_2$} (B3)
      (B3) edge[susp,densely dashed,gray] node[inner sep=1pt,fill=white]{$\beta_3$} (B1)
      (B1) edge[densely dashed,gray] node[inner sep=1pt,fill=white]{$\beta_1$} (B2)
      (B2) edge node[inner sep=1pt,fill=white]{$\tau_2$} (T)
      (T) edge[susp] node[inner sep=2pt,fill=white]{$\tau_3$} (A2)
      (A3) edge[susp] node[inner sep=1pt,fill=white]{$\alpha_3$} (A1)
      (A2) edge node[inner sep=2pt,fill=white]{$\alpha_2$} (A3)
      (A2) edge[densely dashed,gray] node[pos=.42,inner sep=1pt,fill=white]{$\varphi_2$} (B2)
      (A1) edge[densely dashed,gray] node[inner sep=1pt,fill=white]{$\alpha_1$} (A2)

      (S) edge[-,line width=4,white] (B3)
      (S) edge[blue] node[pos=.3,inner sep=1pt,fill=white]{$\nu_2$} (B3)
      (T') edge[-,line width=4,white] (S)
      (T') edge[blue] node[inner sep=1pt,fill=white]{$\nu_1$} (S)
      (S) edge[-,line width=4,white] (T)
      (S) edge[blue] node[inner sep=1pt,fill=white]{$\mu_2$} (T)
      (T) edge[-,line width=4,white] (A3)
      (T) edge[susp,blue] node[pos=.45,inner sep=1pt,fill=white]{$\mu_3$} (A3)
      (A3) edge[-,line width=4,white] (S)
      (A3) edge[blue] node[pos=.55,inner sep=1pt,fill=white]{$\mu_1$} (S);
\end{tikzpicture}
  \caption{Two octahedra glued together by a face. In blue: the ``butterfly.''}
  \label{fig:octahedra}
\end{figure}

\begin{definition}\label{Verdiergood}
	A map of distinguished triangles (\ref{morphism}) is \emph{Verdier good} if $\varphi_3$ factors through an object $S$ as $\nu_2\mu_1$ where $\mu_1$ arises from the octahedron

	\begin{equation}\label{diag:T}
		\begin{tikzpicture}[baseline=(current bounding box.center)]
			\diagram{d}{2em}{3em}{
			  A_1 & A_2 & A_3 & \Sigma A_1\\
			  A_1 & B_2 & S & \Sigma A_1\\
			  A_2 & B_2 & T & \Sigma A_2 \\
			      &     & \Sigma A_3\\
			};
	
			\path[->, font = \scriptsize, auto]
			  (d-1-1) edge node{$\alpha_1$} (d-1-2)
			  (d-1-2) edge node{$\alpha_2$} (d-1-3)
			  (d-1-3) edge node{$\alpha_3$} (d-1-4)
			  
			  (d-2-1) edge node{$\varphi_2\alpha_1$} (d-2-2)
			  (d-2-2) edge node{$\sigma_2$} (d-2-3)
			  (d-2-3) edge node{$\sigma_3$} (d-2-4)
			  
			  (d-3-1) edge node{$\varphi_2$} (d-3-2)
			  (d-3-2) edge node{$\tau_2$} (d-3-3)
			  (d-3-3) edge node{$\tau_3$} (d-3-4)
			  
			  ([xshift=-0.25mm] d-1-1.south) edge[-] ([xshift=-0.25mm]d-2-1.north)
			  ([xshift=0.25mm] d-1-1.south) edge[-] ([xshift=0.25mm]d-2-1.north)
			  (d-1-2) edge node{$\varphi_2$} (d-2-2)
			  (d-1-3) edge[densely dashed] node{$\mu_1$} (d-2-3)
			  ([xshift=-0.25mm] d-1-4.south) edge[-] ([xshift=-0.25mm]d-2-4.north)
			  ([xshift=0.25mm] d-1-4.south) edge[-] ([xshift=0.25mm]d-2-4.north)

			  (d-2-1) edge node{$\alpha_1$} (d-3-1)
			  ([xshift=-0.25mm] d-2-2.south) edge[-] ([xshift=-0.25mm]d-3-2.north)
			  ([xshift=0.25mm] d-2-2.south) edge[-] ([xshift=0.25mm]d-3-2.north)
			  (d-2-3) edge[densely dashed] node{$\mu_2$} (d-3-3)
			  (d-2-4) edge node{$\Sigma\alpha_1$} (d-3-4)

			  (d-3-3) edge[densely dashed] node{$\mu_3$} (d-4-3)
			  (d-3-4) edge node{$\Sigma\alpha_2$} (d-4-3)
			  ;
	\end{tikzpicture}
\end{equation}

	\noindent and $\nu_2$ arises from the octahedron
	
	\begin{equation}\label{diag:T'}
		\begin{tikzpicture}[baseline=(current bounding box.center)]
			\diagram{d}{2em}{3em}{
			  A_1 & B_1 & R & \Sigma A_1\\
			  A_1 & B_2 & S & \Sigma A_1\\
			  B_1 & B_2 & B_3 & \Sigma B_1\\
			      &     & \Sigma R\\ 
			};
	
			\path[->, font = \scriptsize, auto]
			  (d-1-1) edge node{$\varphi_1$} (d-1-2)
			  (d-1-2) edge node{$\rho_2$} (d-1-3)
			  (d-1-3) edge node{$\rho_3$} (d-1-4)
			  
			  (d-2-1) edge node{$\beta_1\varphi_1$} (d-2-2)
			  (d-2-2) edge node{$\sigma_2$} (d-2-3)
			  (d-2-3) edge node{$\sigma_3$} (d-2-4)
			  
			  (d-3-1) edge node{$\beta_1$} (d-3-2)
			  (d-3-2) edge node{$\beta_2$} (d-3-3)
			  (d-3-3) edge node{$\beta_3$} (d-3-4)
			  
			  ([xshift=-0.25mm] d-1-1.south) edge[-] ([xshift=-0.25mm]d-2-1.north)
			  ([xshift=0.25mm] d-1-1.south) edge[-] ([xshift=0.25mm]d-2-1.north)
			  (d-1-2) edge node{$\beta_1$} (d-2-2)
			  (d-1-3) edge[densely dashed] node{$\nu_1$} (d-2-3)
			  ([xshift=-0.25mm] d-1-4.south) edge[-] ([xshift=-0.25mm]d-2-4.north)
			  ([xshift=0.25mm] d-1-4.south) edge[-] ([xshift=0.25mm]d-2-4.north)

			  (d-2-1) edge node{$\varphi_1$} (d-3-1)
			  ([xshift=-0.25mm] d-2-2.south) edge[-] ([xshift=-0.25mm]d-3-2.north)
			  ([xshift=0.25mm] d-2-2.south) edge[-] ([xshift=0.25mm]d-3-2.north)
			  (d-2-3) edge[densely dashed] node{$\nu_2$} (d-3-3)
			  (d-2-4) edge node{$\Sigma\varphi_1$} (d-3-4)

			  (d-3-3) edge[densely dashed] node{$\nu_3$} (d-4-3)
			  (d-3-4) edge node{$\Sigma\rho_2$} (d-4-3)
			  ;
		\end{tikzpicture}
	\end{equation}
\noindent	A fill-in that factors in this way is called a \emph{Verdier good fill-in}.
\end{definition}
Note that $\beta_1\varphi_1=\varphi_2\alpha_1$, so we can, in fact, use the same object $S$ in both octahedra. Since the two octahedral diagrams share the middle row, pictorially the two octahedra will share a face and the two arising distinguished triangles form a \dq{butterfly} with $S$ at its center. Thus $\varphi_3$ should arise as the composition of maps going from the lower lefthand corner to the upper righthand corner in Figure \ref{fig:octahedra}.
Conversely, extending $\varphi_1$ and $\varphi_2$ to distinguished triangles with cones $R$ and $T$ respectively, the other composition through $S$, namely $\mu_2\nu_1$, is a Verdier good fill-in for a morphism between these two new distinguished triangles.

\begin{remark}\label{relationshipsgoodness}
	Given the solid part of the diagram (\ref{diag:fill-in}), grafting together (\ref{diag:T}) and (\ref{diag:T'}) via their shared row and composing the result produces a Verdier good morphism from (TR4) and (TR1)(c). Thus a Verdier good fill-in always exists.

This notion of goodness was originally applied by Verdier in his proof of the so-called $(4\times 4)$-lemma which is described in \cite{Beilinson2018} in order to prove that a middling good fill-in always exists, and said proof entails that Verdier good morphisms are always middling good. Furthermore, it is proven in \cite[Theorem 2.3]{neeman91} that good morphisms are always middling good, so middling goodness is the weakest goodness constraint we can apply to a morphism. It is also known that middling good morphisms are not necessarily good nor Verdier good, however, we can say for certain when a middling good morphism is Verdier good, as Verdier goodness is equivalent to being middling good and satisfying a specific extra condition as stated in \cite[Lemma 3.4]{frankland20}. 
\end{remark}

\section{Higher triangulated categories}

In this section we introduce the generalised notions of Verdier good and middling good morphisms in $n$-angulated categories.

\begin{definition}
	In an $n$-angulated category, a morphism of $n$-angles
	
	\begin{equation}\label{hmorphism}
    \begin{tikzpicture}[baseline=(current bounding box.center)]
      \diagram{d}{2em}{3em}{
        A_1 & A_2 & A_3 & \cdots & A_n & \Sigma A_1\\
        B_1 & B_2 & B_3 & \cdots & B_n & \Sigma B_1\\ 
      };
  
      \path[->, font = \scriptsize, auto]
        (d-1-1) edge node{$\alpha_1$} (d-1-2)
        (d-1-2) edge node{$\alpha_2$} (d-1-3)
        (d-1-3) edge node{$\alpha_3$} (d-1-4)
        (d-1-4) edge node{$\alpha_{n-1}$} (d-1-5)
        (d-1-5) edge node{$\alpha_n$} (d-1-6)
        
        (d-2-1) edge node{$\beta_1$} (d-2-2)
        (d-2-2) edge node{$\beta_2$} (d-2-3)
        (d-2-3) edge node{$\beta_3$} (d-2-4)
        (d-2-4) edge node{$\beta_{n-1}$} (d-2-5)
        (d-2-5) edge node{$\beta_n$} (d-2-6)
        
        (d-1-1) edge node{$\varphi_1$} (d-2-1)
        (d-1-2) edge node{$\varphi_2$} (d-2-2)
        (d-1-3) edge node{$\varphi_3$} (d-2-3)
        (d-1-5) edge node{$\varphi_n$} (d-2-5)
        (d-1-6) edge node{$\Sigma\varphi_1$} (d-2-6)
        ;
    \end{tikzpicture}
  \end{equation}

\noindent	is \emph{Verdier good} if for $3\<i\<n$ the morphism $\varphi_i$ factors as $\nu_{2i}\mu_{1i}$ through an object $S_i$, where the $\mu_{1i}$ arise from an octahedron
	
	\begin{equation}
    \begin{tikzpicture}[baseline=(current bounding box.center)]
      \diagram{d}{2em}{2.5em}{
        A_1 & A_2 & A_3 & A_4 & \cdots & A_{n-1} & A_n & \Sigma A_1\\
        A_1 & B_2 & S_3 & S_4 & \cdots & S_{n-1} & S_n & \Sigma A_1\\ 
        A_2 & B_2 & T_3 & T_4 & \cdots & T_{n-1} & T_n & \Sigma T_2\\
      };
  
      \path[->, font = \scriptsize, auto]
        (d-1-4) edge[densely dashed,out=260,in=30] node[above left,xshift=2pt,pos=.2]{$\lambda_4$} (d-3-3)
        (d-1-7) edge[densely dashed,out=260,in=30] node[above left,xshift=2pt,pos=.2]{$\lambda_n$} (d-3-6)

        (d-1-1) edge node{$\alpha_1$} (d-1-2)
        (d-1-2) edge node{$\alpha_2$} (d-1-3)
        (d-1-3) edge node{$\alpha_3$} (d-1-4)
        (d-1-4) edge node{$\alpha_4$} (d-1-5)
        (d-1-5) edge node{$\alpha_{n-2}$} (d-1-6)
        (d-1-6) edge node{$\alpha_{n-1}$} (d-1-7)
        (d-1-7) edge node{$\alpha_n$} (d-1-8)
        
        (d-2-1) edge node{$\varphi_2\alpha_1$} (d-2-2)
        (d-2-2) edge node{$\sigma_2$} (d-2-3)
        (d-2-3) edge[line width=4,white] (d-2-4)
        (d-2-3) edge node{$\sigma_3$} (d-2-4)
        (d-2-4) edge node{$\sigma_4$} (d-2-5)
        (d-2-5) edge node{$\sigma_{n-2}$} (d-2-6)
        (d-2-6) edge[line width=4,white] (d-2-7)
        (d-2-6) edge node{$\sigma_{n-1}$} (d-2-7)
        (d-2-7) edge node{$\sigma_n$} (d-2-8)

        (d-3-1) edge node{$\varphi_2$} (d-3-2)
        (d-3-2) edge node{$\tau_2$} (d-3-3)
        (d-3-3) edge node{$\tau_3$} (d-3-4)
        (d-3-4) edge node{$\tau_4$} (d-3-5)
        (d-3-5) edge node{$\tau_{n-2}$} (d-3-6)
        (d-3-6) edge node{$\tau_{n-1}$} (d-3-7)
        (d-3-7) edge node{$\tau_n$} (d-3-8)

        ([xshift=-0.25mm] d-1-1.south) edge[-] ([xshift=-0.25mm]d-2-1.north)
        ([xshift=0.25mm] d-1-1.south) edge[-] ([xshift=0.25mm]d-2-1.north)
        (d-1-2) edge node{$\varphi_2$} (d-2-2)
        (d-1-3) edge[densely dashed] node{$\mu_{13}$} (d-2-3)
        (d-1-4) edge[densely dashed] node{$\mu_{14}$} (d-2-4)
        (d-1-6) edge[densely dashed] node{$\mu_{1,n-1}$} (d-2-6)
        (d-1-7) edge[densely dashed] node{$\mu_{1n}$} (d-2-7)
        ([xshift=-0.25mm] d-1-8.south) edge[-] ([xshift=-0.25mm]d-2-8.north)
        ([xshift=0.25mm] d-1-8.south) edge[-] ([xshift=0.25mm]d-2-8.north)
        
        (d-2-1) edge node{$\alpha_1$} (d-3-1)
        ([xshift=-0.25mm] d-2-2.south) edge[-] ([xshift=-0.25mm]d-3-2.north)
        ([xshift=0.25mm] d-2-2.south) edge[-] ([xshift=0.25mm]d-3-2.north)
        (d-2-3) edge[densely dashed] node{$\mu_{23}$} (d-3-3)
        (d-2-4) edge[densely dashed] node{$\mu_{24}$} (d-3-4)
        (d-2-6) edge[densely dashed] node{$\mu_{2,n-1}$} (d-3-6)
        (d-2-7) edge[densely dashed] node{$\mu_{2n}$} (d-3-7)
        (d-2-8) edge node{$\Sigma\alpha_1$} (d-3-8)
        ;
    \end{tikzpicture}
  \end{equation}

	\noindent and the $\nu_{2i}$ arise from an octahedron
	
	\begin{equation}
    \begin{tikzpicture}[baseline=(current bounding box.center)]
      \diagram{d}{2em}{2.5em}{
        A_1 & B_1 & R_3 & R_4 & \cdots & R_{n-1} & R_n & \Sigma A_1\\
        A_1 & B_2 & S_3 & S_4 & \cdots & S_{n-1} & S_n & \Sigma A_1\\ 
        B_1 & B_2 & B_3 & B_4 & \cdots & B_{n-1} & B_n & \Sigma B_2\\
      };
  
      \path[->, font = \scriptsize, auto]
        (d-1-4) edge[densely dashed,out=260,in=30] node[above left,xshift=2pt,pos=.2]{$\gamma_4$} (d-3-3)
        (d-1-7) edge[densely dashed,out=260,in=30] node[above left,xshift=2pt,pos=.2]{$\gamma_n$} (d-3-6)

        (d-1-1) edge node{$\varphi_1$} (d-1-2)
        (d-1-2) edge node{$\rho_2$} (d-1-3)
        (d-1-3) edge node{$\rho_3$} (d-1-4)
        (d-1-4) edge node{$\rho_4$} (d-1-5)
        (d-1-5) edge node{$\rho_{n-2}$} (d-1-6)
        (d-1-6) edge node{$\rho_{n-1}$} (d-1-7)
        (d-1-7) edge node{$\rho_n$} (d-1-8)
        
        (d-2-1) edge node{$\beta_1\varphi_1$} (d-2-2)
        (d-2-2) edge node{$\sigma_2$} (d-2-3)
        (d-2-3) edge[line width=4,white] (d-2-4)
        (d-2-3) edge node{$\sigma_3$} (d-2-4)
        (d-2-4) edge node{$\sigma_4$} (d-2-5)
        (d-2-5) edge node{$\sigma_{n-2}$} (d-2-6)
        (d-2-6) edge[line width=4,white] (d-2-7)
        (d-2-6) edge node{$\sigma_{n-1}$} (d-2-7)
        (d-2-7) edge node{$\sigma_n$} (d-2-8)

        (d-3-1) edge node{$\beta_1$} (d-3-2)
        (d-3-2) edge node{$\beta_2$} (d-3-3)
        (d-3-3) edge node{$\beta_3$} (d-3-4)
        (d-3-4) edge node{$\beta_4$} (d-3-5)
        (d-3-5) edge node{$\beta_{n-2}$} (d-3-6)
        (d-3-6) edge node{$\beta_{n-1}$} (d-3-7)
        (d-3-7) edge node{$\beta_n$} (d-3-8)

        ([xshift=-0.25mm] d-1-1.south) edge[-] ([xshift=-0.25mm]d-2-1.north)
        ([xshift=0.25mm] d-1-1.south) edge[-] ([xshift=0.25mm]d-2-1.north)
        (d-1-2) edge node{$\beta_1$} (d-2-2)
        (d-1-3) edge[densely dashed] node{$\nu_{13}$} (d-2-3)
        (d-1-4) edge[densely dashed] node{$\nu_{14}$} (d-2-4)
        (d-1-6) edge[densely dashed] node{$\nu_{1,n-1}$} (d-2-6)
        (d-1-7) edge[densely dashed] node{$\nu_{1n}$} (d-2-7)
        ([xshift=-0.25mm] d-1-8.south) edge[-] ([xshift=-0.25mm]d-2-8.north)
        ([xshift=0.25mm] d-1-8.south) edge[-] ([xshift=0.25mm]d-2-8.north)
        
        (d-2-1) edge node{$\varphi_1$} (d-3-1)
        ([xshift=-0.25mm] d-2-2.south) edge[-] ([xshift=-0.25mm]d-3-2.north)
        ([xshift=0.25mm] d-2-2.south) edge[-] ([xshift=0.25mm]d-3-2.north)
        (d-2-3) edge[densely dashed] node{$\nu_{23}$} (d-3-3)
        (d-2-4) edge[densely dashed] node{$\nu_{24}$} (d-3-4)
        (d-2-6) edge[densely dashed] node{$\nu_{2,n-1}$} (d-3-6)
        (d-2-7) edge[densely dashed] node{$\nu_{2n}$} (d-3-7)
        (d-2-8) edge node{$\Sigma\varphi_1$} (d-3-8)
        ;
    \end{tikzpicture}
  \end{equation}

	\end{definition}
\begin{remark}
	A Verdier good fill-in always exists by a similar argument to that of triangulated categories. See for example the proof of \cite[Theorem 3.1]{arentz} wherein the morphism constructed is Verdier good.
\end{remark}

\begin{definition}
	A morphism of $n$-angles (\ref{hmorphism}) in an $n$-angulated category, where we label 
	\[
		A_{1i}\coloneqq A_i,\quad A_{2i}\coloneqq B_i,\quad \alpha_{1i}\coloneqq\alpha_i,\quad \alpha_{2i}\coloneqq\beta_i,\quad \mt{and}\quad \varphi_{1i}\coloneqq \varphi_i,
	\]
	is \emph{middling good} if it can be extended to a diagram
	
	\begin{equation}
    \begin{tikzpicture}[baseline=(current bounding box.center)]
      \diagram{d}{2em}{3em}{
        A_{11} & A_{12} & \cdots & A_{1n} & \Sigma A_{11}\\
        A_{21} & A_{22} & \cdots & A_{2n} & \Sigma A_{21}\\
        \vdots & \vdots &  & \vdots & \vdots \\ 
        A_{n1} & A_{n2} & \cdots & A_{nn} & \Sigma A_{n1}\\
        \Sigma A_{11} & \Sigma A_{12} & \cdots & \Sigma A_{1n} & \Sigma^2 A_{11}\\
      };
  
      \path[->, font = \scriptsize, auto]
        (d-1-1) edge node{$\alpha_{11}$} (d-1-2)
        (d-1-2) edge node{$\alpha_{12}$} (d-1-3)
        (d-1-3) edge node{$\alpha_{1,n-1}$} (d-1-4)
        (d-1-4) edge node{$\alpha_{1n}$} (d-1-5)
        
        (d-2-1) edge node{$\alpha_{21}$} (d-2-2)
        (d-2-2) edge node{$\alpha_{22}$} (d-2-3)
        (d-2-3) edge node{$\alpha_{2,n-1}$} (d-2-4)
        (d-2-4) edge node{$\alpha_{2n}$} (d-2-5)
        
        (d-4-1) edge node{$\alpha_{n1}$} (d-4-2)
        (d-4-2) edge node{$\alpha_{n2}$} (d-4-3)
        (d-4-3) edge node{$\alpha_{n,n-1}$} (d-4-4)
        (d-4-4) edge node{$\alpha_{nn}$} (d-4-5)
        
        (d-5-1) edge node{$\Sigma\alpha_{11}$} (d-5-2)
        (d-5-2) edge node{$\Sigma\alpha_{12}$} (d-5-3)
        (d-5-3) edge node{$\Sigma\alpha_{1,n-1}$} (d-5-4)
        (d-5-4) edge node{$\Sigma\alpha_{1n}$} (d-5-5)
        
        (d-1-1) edge node{$\varphi_{11}$} (d-2-1)
        (d-1-2) edge node{$\varphi_{12}$} (d-2-2)
        (d-1-4) edge node{$\varphi_{1n}$} (d-2-4)
        (d-1-5) edge node{$\Sigma\varphi_{11}$} (d-2-5)
        
        (d-2-1) edge node{$\varphi_{21}$} (d-3-1)
        (d-2-2) edge node{$\varphi_{22}$} (d-3-2)
        (d-2-4) edge node{$\varphi_{2n}$} (d-3-4)
        (d-2-5) edge node{$\Sigma\varphi_{21}$} (d-3-5)

        (d-3-1) edge node{$\varphi_{n-1,1}$} (d-4-1)
        (d-3-2) edge node{$\varphi_{n-1,2}$} (d-4-2)
        (d-3-4) edge node{$\varphi_{n-1,n}$} (d-4-4)
        (d-3-5) edge node{$\Sigma\varphi_{n-1,1}$} (d-4-5)

        (d-4-1) edge node{$\varphi_{n1}$} (d-5-1)
        (d-4-2) edge node{$\varphi_{n2}$} (d-5-2)
        (d-4-4) edge node{$\varphi_{nn}$} (d-5-4)
        (d-4-5) edge node{$\Sigma\varphi_{n1}$} (d-5-5)
        ;
        \path (d-4-4) -- (d-5-5) node[midway,xshift=.5em] (b) {$\Omega$};
    \end{tikzpicture}
  \end{equation}
	
	\noindent where the first $n$ columns and $n$ rows are all $n$-angles and in which all squares commute when $n$ is even, and the square $\Omega$ anticommutes when $n$ is odd.
\end{definition}

We should like to observe the relationships of Remark \ref{relationshipsgoodness} in the $n$-angulated setting and, while each goodness condition has a quite clear candidate in this setting as seen above, the arguments and methods used by Verdier, Neeman, and Christensen--Frankland do not easily generalise. Even the proof that Verdier good morphisms are middling good in a 4-angulated category eludes us; a $(5\times 5)$-diagram consisting mostly of $4$-angles is obtainable but, so far, we have been unable to prove whether or not the square of cone objects in the diagram commutes. For now, let us move away from the model free setting in order to say something concrete.


\chapter{Concrete higher triangulations}\label{clustertilting}

It was originally proved in \cite{gko} that an $(n-2)$-cluster tilting subcategory in a triangulated category with suspension $\Sigma$ has an $n$-angulation if it is stable with respect to $\Sigma^{n-2}$. We show in Section 4.1 that not only are Verdier good morphisms also middling good in this class of $n$-angulated categories, in fact, all morphisms of $n$-angles are middling good. In Section 4.2 we discuss exotic $n$-angulated categories.

\section{Cluster tilting categories}

We begin by recalling the relevant definitions to define the $n$-angulated structure on $(n-2)$-cluster tilting categories. Following that, we
 make some specific observations on the nature of $n$-angles in this setting before moving on to the proof of the claim.
\begin{definition}
For a category \cat T and a full subcategory \cat C, an object $X\in\cat T$ has a \emph{left \cat C-approximation} if there exists a morphism $f\in \cat T(X,C)$ with $C\in\cat C$ such that for any object $C'\in \cat C$ and map $g\in \cat T(X,C')$ there exists a map $h\in \cat C(C,C')$ such that $hf=g$.
\vspace{-1em}
\begin{center}
		\begin{tikzpicture}	
		\diagram{d}{2em}{4em}{
		  X & C\\
		  & C'\\ 
		};

		\path[->, font = \scriptsize, auto]
		  (d-1-1) edge node{$f$} (d-1-2)
		  (d-1-1) edge node[below left]{$g$} (d-2-2)
		  (d-1-2) edge[densely dashed] node{$h$} (d-2-2)
		  ;
	\end{tikzpicture}
\end{center}
\vspace{-1em}
\noindent If all objects in \cat T have this property, it is called \emph{covariantly finite}. If $\cat T^{\mt{op}}$ is covariantly finite, then \cat T is \emph{contravariantly finite}. Finally, if \cat T is both covariantly and contravariantly finite, we say that \cat T is \emph{functorially finite.}
\end{definition}

\begin{definition}
Let \cat T be a triangulated category with a class of distinguished triangles $\Delta$ and with an autoequivalence $\Sigma$. A full subcategory $\cat C\subset\cat T$ is \emph{$d$-cluster tilting} if it is functorially finite and it satisfies the condition
\begin{align*}
	\cat C&=\Set[\big]{X\in\cat T\given \cat T(\cat C,\Sigma^iX)=0\mt{ for all }i\in\{1,\dots,d-1\}}\\
		&=\Set[\big]{X\in\cat T\given \cat T(X,\Sigma^i\cat C)=0\mt{ for all }i\in\{1,\dots,d-1\}}.
\end{align*}
\end{definition}	

\noindent Let us now recall the theorem of \cite{gko}.
\begin{theorem}[Geiss--Keller--Oppermann]\label{cluster_is_nang}
Let $(\cat T,\Sigma,\Delta)$ be a triangulated category with an $(n-2)$-cluster tilting subcategory $\cat C$ such that the objects of $\Sigma^{n-2}\cat C$ are contained in \cat C. Then $(\cat C,\Sigma^{n-2},\mathscr N)$ is an $n$-angulated category, where $\mathscr N$ is the class of all sequences
\begin{center}
		\begin{tikzpicture}	
		\diagram{d}{2em}{3em}{
		  A_1 & A_2 & \cdots & A_n & \Sigma^{n-2}A_1\\ 
		};

		\path[->, font = \scriptsize, auto]
		  (d-1-1) edge node{$\alpha_1$} (d-1-2)
		  (d-1-2) edge node{$\alpha_2$} (d-1-3)
		  (d-1-3) edge node{$\alpha_{n-1}$} (d-1-4)
		  (d-1-4) edge node{$\alpha_n$} (d-1-5)
		  ;
	\end{tikzpicture}
\end{center}
such that there exists a diagram
	\begin{equation*}\label{fig:trapez}
	\begin{tikzpicture}[baseline=(current bounding box.center)]	
		\diagram{d}{2.5em}{.7em}{
		   	   & A_2    	 && A_3  		&& A_4 		  && \cdots 		  && A_{n-1}\\
		   A_1 & 	& A_{2.5} &    & A_{3.5} &	  & \cdots &	   & A_{n-1.5} &		& A_n,\\
		   };

		\path[->, font = \scriptsize, auto,inner sep=1pt]
		  (d-1-2) edge node{$\alpha_2$} (d-1-4)
		  (d-1-4) edge node{$\alpha_3$} (d-1-6)
		  (d-1-6) edge node{$\alpha_4$} (d-1-8)
		  (d-1-8) edge node{$\alpha_{n-2}$} (d-1-10)

		  (d-2-1) edge node{$\alpha_1$} (d-1-2)
		  (d-1-2) edge node{$\alpha_2'$} (d-2-3)
		  (d-2-3) edge[susp] node{$\alpha^2_n$} (d-2-1)

		  (d-2-3) edge node{$\alpha_2''$} (d-1-4)
		  (d-1-4) edge node{$\alpha_3'$} (d-2-5)
		  (d-2-5) edge[susp] node{$\alpha^3_n$} (d-2-3)

		  (d-2-5) edge node{$\alpha_3''$} (d-1-6)
		  (d-2-7) edge[susp] node{$\alpha_n^4$} (d-2-5)

		  (d-2-9) edge[susp] node{$\alpha_n^{n-2}$} (d-2-7)
		  
		  (d-2-9) edge node{$\alpha_{n-2}''$} (d-1-10)
		  (d-1-10) edge node{$\alpha_{n-1}$} (d-2-11)
		  (d-2-11) edge[susp] node{$\alpha_n^{n-1}$} (d-2-9)
		  ;
	\end{tikzpicture}
	\end{equation*}\\
\noindent which we will refer to as a \emph{the trapez diagram of an $n$-angle}, with $A_i\in\cat T$ for $i\not\in\Z$, such that 
\[
	(\alpha_i'',\alpha_{i+1}',\alpha_n^{i+1})\in\Delta\quad \mt{\itshape for}\quad i\in\{1,\dots,n-2\},
\]
  where $\alpha_1''\coloneqq\alpha_1$ and $\alpha_{n-1}'\coloneqq\alpha_{n-1}$, and such that 
 \[
 \alpha_i=\alpha_i''\alpha_i'\quad\mt{\itshape for}\quad i\in\{2,\dots,n-2\}
 \]
 and
 \[
\alpha_n=(\Sigma^{n-3}\alpha_n^2)(\Sigma^{n-4}\alpha_n^3)\cdots(\Sigma\alpha_n^{n-2})\alpha_n^{n-1}.
\]

We call the objects $A_{i+1.5}$ the \emph{supporting objects of the $n$-angle} for $i\in\{1,\dots,n-3\}$.
\end{theorem}
\begin{remark}\label{doubledashisapprox}
Assuming $n>3$, in proving that the category \cat C in Theorem \ref{cluster_is_nang} satisfies axiom (N1)(c), one can construct a \cat C-resolution of $A_1$ inductively, starting with $\alpha_1$, by applying (TR1)(c) for $\alpha_i''$ ($\alpha_1''=\alpha_1$) and then taking the left \cat C-approximations of the arising cone of $\alpha_i''$, i.e. the supporting objects, to obtain the maps $\alpha_i'$. Proving that this process terminates, i.e.\ that the cone of $\alpha_{n-2}''$ is an object in \cat C, is, in our opinion, not trivial but we will not reproduce the argument here. Conversely, given a trapez diagram as in Theorem \ref{cluster_is_nang}, it turns out that $\alpha_{i+1}''\colon A_{i+1.5}\to A_{i+2}$ must, in fact, be a left approximation. To see this, fix an $i$ such that $1\<i\<n-3$. First, given a morphism $\gamma\colon A_{i+1.5}\to C$ where $C\in\cat C$, since 
\[
	\gamma\alpha_{i+1}'\alpha_i=\gamma(\alpha_{i+1}'\alpha_i'')\alpha_i'=0,
\]
we obtain the solid diagram
\begin{center}
		\begin{tikzpicture}	
		\diagram{d}{2em}{3em}{
		  A_{i} & A_{i+1} & A_{i+2} & \cdots & \Sigma A_{i-1} & \Sigma A_{i}\\
		  0 & C & C & \cdots & 0 & 0\\ 
		};

		\path[->, font = \scriptsize, auto]
		  (d-1-1) edge node{$\alpha_{i}$} (d-1-2)
		  (d-1-2) edge node{$\alpha_{i+1}$} (d-1-3)
		  (d-1-3) edge node{$\alpha_{i+2}$} (d-1-4)
		  (d-1-4) edge node{$(-1)^n\Sigma\alpha_{i-1}$} (d-1-5)
		  (d-1-5) edge node{$(-1)^n\Sigma\alpha_i$} (d-1-6)
		  
		  (d-2-1) edge (d-2-2)
		  ([yshift=-0.25mm] d-2-2.east) edge[-] ([yshift=-0.25mm]d-2-3.west)
		  ([yshift=0.25mm] d-2-2.east) edge[-] ([yshift=0.25mm]d-2-3.west)
		  (d-2-3) edge (d-2-4)
		  (d-2-4) edge (d-2-5)
		  (d-2-5) edge (d-2-6)

		  (d-1-1) edge (d-2-1)
		  (d-1-2) edge node{$\gamma\alpha_{i+1}'$} (d-2-2)
		  (d-1-3) edge[densely dashed] node{$\bar\gamma$} (d-2-3)
		  (d-1-5) edge[densely dashed] (d-2-5)
		  (d-1-6) edge (d-2-6)
		  ;
	\end{tikzpicture}
\end{center}
and applying (N3) we gain a dashed map $\bar\gamma\in \cat C(A_{i+2}, C)$ such that
\[
		\bar\gamma\alpha_{i+1}=\gamma\alpha_{i+1}'.
\]
From this identity, as $\alpha_{i+1}=\alpha_{i+1}''\alpha_{i+1}'$, we can conclude that 
\[
	(\bar\gamma\alpha_{i+1}''-\gamma)\alpha_{i+1}'=0,
\]
so since $\alpha_n^{i+1}\colon A_{i+1.5}\to \Sigma A_{i+.5}$ is a weak cokernel of $\alpha_{i+1}'$ in \cat T, there is a map
\[
	\Sigma \gamma_1\colon \Sigma A_{i+.5}\to C
\]
satisfying $(\Sigma \gamma_1)\alpha_n^{i+1}=\bar\gamma\alpha_{i+1}''-\gamma$. Recall that \cat C is closed under $\Sigma^{n-2}$, so $\Sigma^{n-2}A_i\in\cat C$, and therefore we get that
\[
	\gamma_1\alpha_i'\in \cat T(A_i,\Sigma^{-1}C)\cong \cat T(\Sigma^{n-2}A_i,\Sigma^{n-3}C)=0
\]
by the cluster tilting condition. Note that the cluster tilting condition only guarantees this because $n-2>1$, and this is why we only consider $(n-2)$-cluster tilting categories for $n>3$.
Since $\alpha^{i}_n$ is a weak cone of $\alpha_i'$, we again obtain a map
\[
	\Sigma \gamma_2\colon \Sigma A_{i-.5}\to \Sigma^{-1}C
\]
such that $(\Sigma \gamma_2)\alpha_n^i=\gamma_1$. We repeat this process until we get a map 
\[
	\Sigma \gamma_i\colon \Sigma A_1\to \Sigma^{1-i}C
\]
such that $(\Sigma \gamma_i)\alpha_n^2=\gamma_{i-1}$. Inductively we see that
\[
	(\Sigma^{i-1}\gamma_i)(\Sigma^{i-2}\alpha_n^2)(\Sigma^{i-3}\alpha_n^3)\cdots(\Sigma\alpha_n^{i-1})\alpha_n^{i}=\gamma_1
\]
but since $1\<i\<n-3$ it follows that $1\<n-i-2\< n-3$, so
\[
	\gamma_i\in \cat T(A_1,\Sigma^{-i}C)\cong \cat T(\Sigma^{n-2}A_1,\Sigma^{n-i-2}C)=0.
\]
The important difference between $\gamma_i$ and $\gamma_k$ for $k<i$ here, is that the source of $\gamma_i$ is an object in \cat C which is not necessarily true of $\gamma_k$. Thus, it follows that $\gamma_1=0$ and we conclude that
\[
	\bar\gamma\alpha_{i+1}''-\gamma=(\Sigma \gamma_1)\alpha_n^{i+1}=0,
\]
i.e.\ for all maps $\gamma\colon A_{i+1.5}\to C$ where $C\in\cat C$, there exists a map $\bar\gamma\colon A_{i+2}\to C$ in \cat C such that $\bar\gamma\alpha_{i+1}''=\gamma$, and therefore $\alpha_{i+1}''\colon A_{i+1.5}\to A_{i+2}$ is a left \cat C-approximation.
\end{remark}

\begin{remark}
Since any $n$-angle in a cluster tilting subcategory is made up of left approximations, we may, similarly to the previous remark, freely construct a morphism of $n$-angles from a square
\begin{center}
 		\begin{tikzpicture}	
 		\diagram{d}{2em}{3em}{
 		  A_1 & A_2\\
 		  B_1 & B_2\\ 
 		};
 
 		\path[->, font = \scriptsize, auto]
 		  (d-1-1) edge node{$\alpha_1$} (d-1-2)
 		  (d-2-1) edge node{$\beta_1$} (d-2-2)
 		  (d-1-1) edge node{$\varphi_1$} (d-2-1)
 		  (d-1-2) edge node{$\varphi_2$} (d-2-2)
 		  ;
 	\end{tikzpicture}
 
 \end{center}
  by iterately applying (TR3) and the universal property of the left approximations to the trapez prism diagram as illustrated in Figure \ref{morphismTrapez1}. The entire diagram of Figure \ref{morphismTrapez1} therefore commutes by construction.
Conversely, given a morphism of $n$-angles, any fill-in between the support objects that arise from (TR3) applied to the constituent distinguished triangles, constitutes a commutative trapez diagram. We state this more precisely as a lemma.
  \begin{figure}[!h]
  \centering
		\begin{tikzpicture}	
		\diagram{d}{2.5em}{.9em}{
		   	   & A_2    	 && A_3  		&& A_4 		  && \cdots 		  && A_{n-1}\\
		   A_1 & 	& A_{2.5} &    & A_{3.5} &	  & \cdots &	   & A_{n-1.5} &		& A_n\\
		   	   & B_2    	 && B_3  		&& B_4 		  && \cdots 		  && B_{n-1}\\
		   B_1 & 	& B_{2.5} &    & B_{3.5} &	  & \cdots &	   & B_{n-1.5} &		& B_n.\\
		};

		\path[->, font = \scriptsize, auto,inner sep=1pt]
		  (d-1-2) edge node{$\alpha_2$} (d-1-4)
		  (d-1-4) edge node{$\alpha_3$} (d-1-6)
		  (d-1-6) edge node{$\alpha_4$} (d-1-8)
		  (d-1-8) edge node{$\alpha_{n-2}$} (d-1-10)

		  (d-3-2) edge node[pos=.25]{$\beta_2$} (d-3-4)
		  (d-3-4) edge node[pos=.25]{$\beta_3$} (d-3-6)
		  (d-3-6) edge node{$\beta_4$} (d-3-8)
		  (d-3-8) edge node[pos=.25]{$\beta_{n-2}$} (d-3-10)

		  (d-1-2) edge node[pos=.28,fill=white,xshift=-.4em]{$\varphi_{2}$} (d-3-2)
		  (d-1-4) edge[densely dashed] node[pos=.28,fill=white,xshift=-.4em]{$\varphi_3$} (d-3-4)
		  (d-1-6) edge[densely dashed] node[pos=.28,fill=white,xshift=-.4em]{$\varphi_4$} (d-3-6)
		  (d-1-10) edge[densely dashed] node[pos=.28,fill=white,xshift=-.6em]{$\varphi_{n-1}$} (d-3-10)

		  (d-2-1) edge node[pos=.4]{$\varphi_{1}$} (d-4-1)
		  (d-2-3) edge[-,white,line width=4] (d-4-3)
		  (d-2-3) edge[densely dotted] node[pos=.25]{$\varphi_{2.5}$} (d-4-3)
		  (d-2-5) edge[-,white,line width=4] (d-4-5)
		  (d-2-5) edge[densely dotted] node[pos=.25]{$\varphi_{3.5}$} (d-4-5)
		  (d-2-9) edge[-,white,line width=4] (d-4-9)
		  (d-2-9) edge[densely dotted] node[pos=.25]{$\varphi_{n-1.5}$} (d-4-9)
		  (d-2-11) edge[densely dotted] node{$\varphi_n$} (d-4-11)

		  (d-2-1) edge node{$\alpha_1$} (d-1-2)
		  (d-1-2) edge node{$\alpha_2'$} (d-2-3)
		  (d-2-3) edge[-,white,line width=4] (d-2-1)
		  (d-2-3) edge[susp] node[pos=.25,yshift=-3pt]{$\alpha^2_n$} (d-2-1)

		  (d-2-3) edge node{$\alpha_2''$} (d-1-4)
		  (d-1-4) edge node{$\alpha_3'$} (d-2-5)
		  (d-2-5) edge[-,white,line width=4] (d-2-3)
		  (d-2-5) edge[susp] node[pos=.25,yshift=-3pt]{$\alpha^3_n$} (d-2-3)

		  (d-2-5) edge node{$\alpha_3''$} (d-1-6)
		  (d-2-7) edge[-,white,line width=4] (d-2-5)
		  (d-2-7) edge[susp] node[pos=.25,yshift=-3pt]{$\alpha_n^4$} (d-2-5)

		  (d-2-9) edge[susp] node[yshift=-3pt]{$\alpha_n^{n-2}$} (d-2-7)
		  
		  (d-2-9) edge node{$\alpha_{n-2}''$} (d-1-10)
		  (d-1-10) edge node{$\alpha_{n-1}$} (d-2-11)
		  (d-2-11) edge[-,white,line width=4] (d-2-9)
		  (d-2-11) edge[susp] node[pos=.25,yshift=-3pt]{$\alpha_n^{n-1}$} (d-2-9)

		  (d-4-1) edge node[pos=.6]{$\beta_1$} (d-3-2)
		  (d-3-2) edge node[pos=.4]{$\beta_2'$} (d-4-3)
		  (d-4-3) edge[-,white,line width=4] (d-4-1)
		  (d-4-3) edge[susp] node{$\beta^2_n$} (d-4-1)

		  (d-4-3) edge node[pos=.6]{$\beta_2''$} (d-3-4)
		  (d-3-4) edge node[pos=.4]{$\beta_3'$} (d-4-5)
		  (d-4-5) edge[-,white,line width=4] (d-4-3)
		  (d-4-5) edge[susp] node{$\beta^3_n$} (d-4-3)

		  (d-4-5) edge node[pos=.6]{$\beta_3''$} (d-3-6)
		  (d-4-7) edge[-,white,line width=4] (d-4-5)
		  (d-4-7) edge[susp] node{$\beta_n^4$} (d-4-5)

		  (d-4-9) edge[susp] node{$\beta_n^{n-2}$} (d-4-7)
		  
		  (d-4-9) edge node[pos=.65]{$\beta_{n-2}''$} (d-3-10)
		  (d-3-10) edge node[pos=.35]{$\beta_{n-1}$} (d-4-11)
		  (d-4-11) edge[-,white,line width=4] (d-4-9)
		  (d-4-11) edge[susp] node{$\beta_n^{n-1}$} (d-4-9)
		  ;
	\end{tikzpicture}
 
\caption{The dotted morphisms come from (TR3) while the dashed morphisms come from the universal property of left approximations.}\label{morphismTrapez1}
\end{figure}
\end{remark}

\begin{lemma}\label{morphismscomefromtri}
	Let $n>3$ be an integer. Suppose we have a morphism of $n$-angles in an $n$-angulated $(n-2)$-cluster tilting subcategory
	  \begin{center}
	  	\begin{tikzpicture}	
		\diagram{d}{2em}{3em}{
		  A_1 & A_2 & A_3 & \cdots & A_n & \Sigma^{n-2} A_1\\
		  B_1 & B_2 & B_3 & \cdots & B_n & \Sigma^{n-2} B_1.\\ 
		};

		\path[->, font = \scriptsize, auto]
		  (d-1-1) edge node{$\alpha_1$} (d-1-2)
		  (d-1-2) edge node{$\alpha_2$} (d-1-3)
		  (d-1-3) edge node{$\alpha_3$} (d-1-4)
		  (d-1-4) edge node{$\alpha_{n-1}$} (d-1-5)
		  (d-1-5) edge node{$\alpha_n$} (d-1-6)
		  
		  (d-2-1) edge node{$\beta_1$} (d-2-2)
		  (d-2-2) edge node{$\beta_2$} (d-2-3)
		  (d-2-3) edge node{$\beta_3$} (d-2-4)
		  (d-2-4) edge node{$\beta_{n-1}$} (d-2-5)
		  (d-2-5) edge node{$\beta_n$} (d-2-6)

		  (d-1-1) edge node{$\varphi_1$} (d-2-1)
		  (d-1-2) edge node{$\varphi_2$} (d-2-2)
		  (d-1-3) edge node{$\varphi_3$} (d-2-3)
		  (d-1-5) edge node{$\varphi_n$} (d-2-5)
		  (d-1-6) edge node{$\Sigma\varphi_1$} (d-2-6)
		  ;
	\end{tikzpicture}
		\end{center}
	where the upper and lower row have support objects $A_{i+1.5}$ for $i\in\{1,\dots,n-3\}$ and $B_{i+1.5}$ for $i\in\{1,\dots,n-3\}$ respectively. If we inductively fill in the diagrams
	\begin{center}
		\begin{tikzpicture}	
			\diagram{d}{2.5em}{3em}{
			  A_{i+1.5} & A_{i+2} & A_{i+2.5} & \Sigma A_{i+1.5}\\ 
			  B_{i+1.5} & B_{i+2} & B_{i+2.5} & \Sigma B_{i+1.5}\\
			};
	
			\path[->, font = \scriptsize, auto]
			  (d-1-1) edge node{$\alpha_{i+1}''$} (d-1-2)
			  (d-1-2) edge node{$\alpha_{i+2}'$} (d-1-3)
			  (d-1-3) edge node{$\alpha_n^{i+2}$} (d-1-4)
			  
			  (d-2-1) edge node{$\beta_{i+1}''$} (d-2-2)
			  (d-2-2) edge node{$\beta_{i+2}'$} (d-2-3)
			  (d-2-3) edge node{$\beta_n^{i+2}$} (d-2-4)

			  (d-1-1) edge node{$\varphi_{i+1.5}$} (d-2-1)
			  (d-1-2) edge node{$\varphi_{i+2}$} (d-2-2)
			  (d-1-3) edge[densely dashed] node{$\varphi_{i+2.5}$} (d-2-3)
			  (d-1-4) edge node{$\Sigma\varphi_{i+1.5}$} (d-2-4)
			  ;
		\end{tikzpicture}
	\end{center}
	starting with $i=0$ and applying the conventions $A_{1.5}\coloneqq A_1$, $B_{1.5}\coloneqq B_1$, $\alpha_1''\coloneqq\alpha_1$, and $\beta_1''\coloneqq\beta_1$, the trapez prism diagram
	
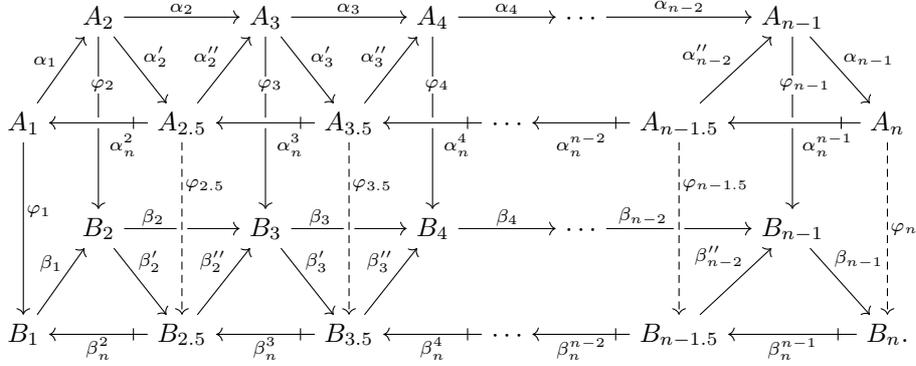
\begin{center}
			\begin{tikzpicture}	
		\diagram{d}{2.5em}{.9em}{
		   	   & A_2    	 && A_3  		&& A_4 		  && \cdots 		  && A_{n-1}\\
		   A_1 & 	& A_{2.5} &    & A_{3.5} &	  & \cdots &	   & A_{n-1.5} &		& A_n\\
		   	   & B_2    	 && B_3  		&& B_4 		  && \cdots 		  && B_{n-1}\\
		   B_1 & 	& B_{2.5} &    & B_{3.5} &	  & \cdots &	   & B_{n-1.5} &		& B_n.\\
		};

		\path[->, font = \scriptsize, auto,inner sep=1pt]
		  (d-1-2) edge node{$\alpha_2$} (d-1-4)
		  (d-1-4) edge node{$\alpha_3$} (d-1-6)
		  (d-1-6) edge node{$\alpha_4$} (d-1-8)
		  (d-1-8) edge node{$\alpha_{n-2}$} (d-1-10)

		  (d-3-2) edge node[pos=.25]{$\beta_2$} (d-3-4)
		  (d-3-4) edge node[pos=.25]{$\beta_3$} (d-3-6)
		  (d-3-6) edge node{$\beta_4$} (d-3-8)
		  (d-3-8) edge node[pos=.25]{$\beta_{n-2}$} (d-3-10)

		  (d-1-2) edge node[pos=.28,fill=white,xshift=-.4em]{$\varphi_{2}$} (d-3-2)
		  (d-1-4) edge node[pos=.28,fill=white,xshift=-.4em]{$\varphi_3$} (d-3-4)
		  (d-1-6) edge node[pos=.28,fill=white,xshift=-.4em]{$\varphi_4$} (d-3-6)
		  (d-1-10) edge node[pos=.28,fill=white,xshift=-.6em]{$\varphi_{n-1}$} (d-3-10)

		  (d-2-1) edge node[pos=.4]{$\varphi_{1}$} (d-4-1)
		  (d-2-3) edge[-,white,line width=4] (d-4-3)
		  (d-2-3) edge[densely dashed] node[pos=.25]{$\varphi_{2.5}$} (d-4-3)
		  (d-2-5) edge[-,white,line width=4] (d-4-5)
		  (d-2-5) edge[densely dashed] node[pos=.25]{$\varphi_{3.5}$} (d-4-5)
		  (d-2-9) edge[-,white,line width=4] (d-4-9)
		  (d-2-9) edge[densely dashed] node[pos=.25]{$\varphi_{n-1.5}$} (d-4-9)
		  (d-2-11) edge[densely dashed] node{$\varphi_n$} (d-4-11)

		  (d-2-1) edge node{$\alpha_1$} (d-1-2)
		  (d-1-2) edge node{$\alpha_2'$} (d-2-3)
		  (d-2-3) edge[-,white,line width=4] (d-2-1)
		  (d-2-3) edge[susp] node[pos=.25,yshift=-3pt]{$\alpha^2_n$} (d-2-1)

		  (d-2-3) edge node{$\alpha_2''$} (d-1-4)
		  (d-1-4) edge node{$\alpha_3'$} (d-2-5)
		  (d-2-5) edge[-,white,line width=4] (d-2-3)
		  (d-2-5) edge[susp] node[pos=.25,yshift=-3pt]{$\alpha^3_n$} (d-2-3)

		  (d-2-5) edge node{$\alpha_3''$} (d-1-6)
		  (d-2-7) edge[-,white,line width=4] (d-2-5)
		  (d-2-7) edge[susp] node[pos=.25,yshift=-3pt]{$\alpha_n^4$} (d-2-5)

		  (d-2-9) edge[susp] node[yshift=-3pt]{$\alpha_n^{n-2}$} (d-2-7)
		  
		  (d-2-9) edge node{$\alpha_{n-2}''$} (d-1-10)
		  (d-1-10) edge node{$\alpha_{n-1}$} (d-2-11)
		  (d-2-11) edge[-,white,line width=4] (d-2-9)
		  (d-2-11) edge[susp] node[pos=.25,yshift=-3pt]{$\alpha_n^{n-1}$} (d-2-9)

		  (d-4-1) edge node[pos=.6]{$\beta_1$} (d-3-2)
		  (d-3-2) edge node[pos=.4]{$\beta_2'$} (d-4-3)
		  (d-4-3) edge[-,white,line width=4] (d-4-1)
		  (d-4-3) edge[susp] node{$\beta^2_n$} (d-4-1)

		  (d-4-3) edge node[pos=.6]{$\beta_2''$} (d-3-4)
		  (d-3-4) edge node[pos=.4]{$\beta_3'$} (d-4-5)
		  (d-4-5) edge[-,white,line width=4] (d-4-3)
		  (d-4-5) edge[susp] node{$\beta^3_n$} (d-4-3)

		  (d-4-5) edge node[pos=.6]{$\beta_3''$} (d-3-6)
		  (d-4-7) edge[-,white,line width=4] (d-4-5)
		  (d-4-7) edge[susp] node{$\beta_n^4$} (d-4-5)

		  (d-4-9) edge[susp] node{$\beta_n^{n-2}$} (d-4-7)
		  
		  (d-4-9) edge node[pos=.65]{$\beta_{n-2}''$} (d-3-10)
		  (d-3-10) edge node[pos=.35]{$\beta_{n-1}$} (d-4-11)
		  (d-4-11) edge[-,white,line width=4] (d-4-9)
		  (d-4-11) edge[susp] node{$\beta_n^{n-1}$} (d-4-9)
		  ;
	\end{tikzpicture}
 
	\end{center}
	will commute.
\end{lemma}

\begin{proof}
We progress in a similar fashion to Remark \ref{doubledashisapprox}. 
Fill in the diagram
\begin{equation}\label{trimorphism}
		\begin{tikzpicture}[baseline=(current bounding box.center)]
			\diagram{d}{2.5em}{3em}{
			  A_{1} & A_{2} & A_{2.5} & \Sigma A_{1}\\ 
			  B_{1} & B_{2} & B_{2.5} & \Sigma B_{1}\\
			};
	
			\path[->, font = \scriptsize, auto]
			  (d-1-1) edge node{$\alpha_{1}$} (d-1-2)
			  (d-1-2) edge node{$\alpha_{2}'$} (d-1-3)
			  (d-1-3) edge node{$\alpha_n^{2}$} (d-1-4)
			  
			  (d-2-1) edge node{$\beta_{1}$} (d-2-2)
			  (d-2-2) edge node{$\beta_{2}'$} (d-2-3)
			  (d-2-3) edge node{$\beta_n^{2}$} (d-2-4)

			  (d-1-1) edge node{$\varphi_{1}$} (d-2-1)
			  (d-1-2) edge node{$\varphi_{2}$} (d-2-2)
			  (d-1-3) edge[densely dashed] node{$\varphi_{2.5}$} (d-2-3)
			  (d-1-4) edge node{$\Sigma\varphi_{1}$} (d-2-4)
			  ;
		\end{tikzpicture}
	\end{equation}
however you like such that the diagram commutes. Of course, there always is at least one choice by (TR3).
Drawing in $\varphi_{2.5}$ in the trapez prism, we know that the solid part commutes along with the two new squares that arise from (\ref{trimorphism}) but, possibly, we might not have the commutativity relation $\varphi_3\alpha_2''=\beta_2''\varphi_{2.5}$. To see that this identity also holds, note that
\[
	\varphi_3\alpha_2''\alpha_2'=\varphi_3\alpha_2=\beta_2\varphi_2=\beta_2''\beta_2'\varphi_2=\beta_2''\varphi_{2.5}\alpha_2'
\]
so
\[
	(\beta_2''\varphi_{2.5}-\varphi_3\alpha_2'')\alpha_2'=0,
\]
and thus there is a morphism
\[
	\Sigma \gamma\colon \Sigma A_1\to B_3
\]
such that $(\Sigma \gamma)\alpha_n^2=\beta_2''\varphi_{2.5}-\varphi_3\alpha_2''$. However, as before, $\Sigma \gamma=0$ since
\[
	\cat T(\Sigma A_1,B_3)\cong \cat T(\Sigma^{n-2}A_1,\Sigma^{n-3}B_3)=0,
\]
 so we conclude that the commutativity relation holds. Again, this follows from the fact that $n-2>1$. Now suppose that $\varphi_{i+1}\alpha_i''=\beta_i''\varphi_{i+.5}$ for some fixed $i\in\{2,\dots,n-3\}$. To see that $\varphi_{i+2}\alpha_{i+1}''=\beta_{i+1}''\varphi_{i+1.5}$, we note again that 
\[
	(\beta_{i+1}''\varphi_{i+1.5}-\varphi_{i+2}\alpha_{i+1}'')\alpha_{i+1}'=0
\]
so there is a map
\[
	\Sigma \gamma_1\colon \Sigma A_{i+.5}\to B_{i+2}
\]
such that $(\Sigma \gamma_1)\alpha_n^{i+1}=\beta_{i+1}''\varphi_{i+1.5}-\varphi_{i+2}\alpha_{i+1}''$, and
\[
	\Sigma (\gamma_1\alpha_i')\in \cat T(\Sigma A_i, B_{i+2})\cong \cat T(\Sigma^{n-2}A_i,\Sigma^{n-3}B_{i+2})=0,
\]
so from the weak cokernel property of $\alpha_n^i$ we get another map
\[
	\Sigma \gamma_2\colon \Sigma A_{i-.5}\to \Sigma^{-1}B_{i+2}
\]
such that $(\Sigma \gamma_2)\alpha_n^i=\gamma_1$. Now repeating this process exactly as we did in Remark \ref{doubledashisapprox} but with $B_{i+2}$ in place of $C$, we conclude that
\[
	0=(\Sigma^{i-1}\gamma_i)(\Sigma^{i-2}\alpha_n^2)(\Sigma^{i-3}\alpha_n^3)\cdots(\Sigma\alpha_n^{i-1})\alpha_n^{i}=\gamma_1
\]
since 
\[
	\gamma_i\in\cat T(A_1,\Sigma^{-i}B_{i+2})\cong \cat T(\Sigma^{n-2}A_1,\Sigma^{n-i-2}B_{i+2})=0
\]
and so
\[
	\varphi_{i+2}\alpha_{i+1}''=\beta_{i+1}''\varphi_{i+.5}
\]
for all $i\in\{2,\dots,n-3\}$.
\end{proof}

We have now arrived at our main result of this paper. Since all $n$-angles in $(n-2)$-cluster tilting subcategories are built from triangles in the ambient triangulated category, diagrams witnessing middling good morphisms can be built up from the triangulated structure.

\begin{theorem}\label{maintheorem}
For $n>3$, all morphisms of $n$-angles in an $(n-2)$-cluster tilting $n$-angulated category are middling good.
\end{theorem}
\paragraph{Notation guide for the proof:} We construct an $(n+1)\times(n+1)$ diagram of objects labeled $A_{ij}$ with $A_{i\bullet}$ and $A_{\bullet j}$ representing $n$-angles and, in particular, when $i$ or $j$ are not integers, $A_{ij}$ represents a support object. We will also adopt the convention that horisontal morphisms are labeled $\alpha_{ij}$ and vertical morphisms are labeled $\varphi_{ij}$ with $ij$ being the index of the source object in both cases. Thus $\alpha_{i\bullet}$ are the structure morphisms of $A_{i\bullet}$ while $\varphi_{\bullet j}$ are the structure morphisms for $A_{\bullet j}$. This will be put aside in the case of support objects, though, as we in this case opt to follow the convention of writing $\alpha_{*,*}''$, $\alpha_{*,*}'$, and $\alpha_{*n}^*$ similarly to what we have been doing so far.
\begin{prf}
We begin with a morphism of $n$-angles, i.e.\ a commutative diagram
\begin{center}
		\begin{tikzpicture}	
		\diagram{d}{2em}{3em}{
		  A_{11} & A_{12} & A_{13} & \cdots & A_{1n} & \Sigma A_{11}\\
		  A_{21} & A_{22} & A_{23} & \cdots & A_{2n} & \Sigma A_{21}\\ 
		};

		\path[->, font = \scriptsize, auto]
		  (d-1-1) edge node{$\alpha_{11}$} (d-1-2)
		  (d-1-2) edge node{$\alpha_{12}$} (d-1-3)
		  (d-1-3) edge node{$\alpha_{13}$} (d-1-4)
		  (d-1-4) edge node{$\alpha_{1,n-1}$} (d-1-5)
		  (d-1-5) edge node{$\alpha_{1n}$} (d-1-6)
		  
		  (d-2-1) edge node{$\alpha_{21}$} (d-2-2)
		  (d-2-2) edge node{$\alpha_{22}$} (d-2-3)
		  (d-2-3) edge node{$\alpha_{23}$} (d-2-4)
		  (d-2-4) edge node{$\alpha_{2,n-1}$} (d-2-5)
		  (d-2-5) edge node{$\alpha_{2n}$} (d-2-6)

		  (d-1-1) edge node{$\varphi_{11}$} (d-2-1)
		  (d-1-2) edge node{$\varphi_{12}$} (d-2-2)
		  (d-1-3) edge node{$\varphi_{13}$} (d-2-3)
		  (d-1-5) edge node{$\varphi_{1n}$} (d-2-5)
		  (d-1-6) edge node{$\Sigma\varphi_{11}$} (d-2-6)
		  ;
	\end{tikzpicture}
 
\end{center}
such that the rows are $n$-angles.
Since we, by Lemma \ref{morphismscomefromtri}, can choose any fill-in morphism we like between the support objects of $A_{1\bullet}$ and $A_{2\bullet}$, we can choose them to be middling good. For short, we may say that the morphism has \emph{middling good support}. For each morphism of triangles, because the third morphism is a middling good fill-in, we have the diagram
\begin{equation}\label{gamma1i}
		\begin{tikzpicture}[baseline=(current bounding box.center)]
		\diagram{d}{3em}{3em}{
		  A_{1,i+.5} & A_{1,i+1} & A_{1,i+1.5} & \Sigma A_{1,i+.5}\\
		  A_{2,i+.5} & A_{2,i+1} & A_{2,i+1.5} & \Sigma A_{2,i+.5}\\
		  A_{2.5,i+.5} & A_{2.5,i+1} & A_{2.5,i+1.5} & \Sigma A_{2.5,i+.5}\\
		  \Sigma A_{1,i+.5} & \Sigma A_{1,i+1} & \Sigma A_{1,i+1.5} & \Sigma^2 A_{1,i+.5}.\\
		};

		\path[->, font = \scriptsize, auto]
		  (d-1-1) edge node{$\alpha_{1i}''$} (d-1-2)
		  (d-1-2) edge node{$\alpha_{1,i+1}'$} (d-1-3)
		  (d-1-3) edge node{$\alpha_{1n}^{i+1}$} (d-1-4)
		  
		  (d-2-1) edge node{$\alpha_{2i}''$} (d-2-2)
		  (d-2-2) edge node{$\alpha_{2,i+1}'$} (d-2-3)
		  (d-2-3) edge node{$\alpha_{2n}^{i+1}$} (d-2-4)
		  
		  (d-3-1) edge[densely dashed] node{$\alpha_{2.5,i}''$} (d-3-2)
		  (d-3-2) edge[densely dashed] node{$\alpha_{2.5,i+1}'$} (d-3-3)
		  (d-3-3) edge[densely dashed] node{$\alpha_{2.5,n}^{i+1}$} (d-3-4)
		  
		  (d-4-1) edge node{$\Sigma\alpha_{1i}''$} (d-4-2)
		  (d-4-2) edge node{$\Sigma\alpha_{1,i+1}'$} (d-4-3)
		  (d-4-3) edge node{$\Sigma\alpha_{1n}^{i+1}$} (d-4-4)
		  
		  (d-1-1) edge node[pos=.4]{$\varphi_{1,i+.5}$} (d-2-1)
		  (d-1-2) edge node[pos=.4]{$\varphi_{1,i+1}$} (d-2-2)
		  (d-1-3) edge node[pos=.4]{$\varphi_{1,i+1.5}$} (d-2-3)
		  (d-1-4) edge node[pos=.4]{$\Sigma\varphi_{1,i+.5}$} (d-2-4)
		  
		  (d-2-1) edge[densely dotted] node[pos=.4]{$\varphi_{2,i+.5}'$} (d-3-1)
		  (d-2-2) edge[densely dotted] node[pos=.4]{$\varphi_{2,i+1}'$} (d-3-2)
		  (d-2-3) edge[densely dotted] node[pos=.4]{$\varphi_{2,i+1.5}'$} (d-3-3)
		  (d-2-4) edge node[pos=.4]{$\Sigma\varphi_{2,i+.5}'$} (d-3-4)
		  
		  (d-3-1) edge[densely dotted] node[pos=.4]{$\varphi_{i+.5,n}^2$} (d-4-1)
		  (d-3-2) edge[densely dotted] node[pos=.4]{$\varphi_{i+1,n}^2$} (d-4-2)
		  (d-3-3) edge[densely dotted] node[pos=.4]{$\varphi_{i+1.5,n}^2$} (d-4-3)
		  (d-3-4) edge node[pos=.4]{$\Sigma\varphi_{i+.5,n}^2$} (d-4-4)
		  ;
		  \path (d-4-3) -- (d-3-4) node[pos=.6,xshift=.5em] (b) {$\Omega_{1i}$};
	\end{tikzpicture}
\end{equation}
The dotted morphisms represent the new distinguished triangles we get from (TR1)(c) and the dashed morphisms represent the new distinguished triangle we get from middling goodness. We now fix $i$. Recall that the square $\Omega_{1i}$ anticommutes and that the morphism $\alpha_{2.5,i}''\colon A_{2.5,i+.5}\to A_{2.5,i+1}$ is a middling good fill-in between the two left-most vertical distinguished triangles in the diagram. By (N1)(c) we can produce an $n$-angle $A_{\bullet,i+.5}$ with base $\varphi_{1,i+.5}\colon A_{1,i+.5}\to A_{2,i+.5}$, and by uniqueness of cones in triangulated categories, we may assume that the object $A_{2.5,i+.5}$ is the initial support object for $A_{\bullet,i+.5}$. Similarly, we construct an $n$-angle $A_{\bullet,i+1}$ with initial support object $A_{2.5,i+1}$. We can further construct a morphism of $n$-angles $A_{\bullet,i+.5}\to A_{\bullet,i+1}$ such that the morphism between the initial support objects is $\alpha_{2.5,i}''$ by repeatedly applying the left approximation property and (TR3), and we may assume that the fill-ins between the remaining support objects are middling good.
Therefore, we in general get middling good diagrams
\begin{figure}[!h]
  \centering
		\begin{tikzpicture}	
		\diagram{d}{3em}{3.5em}{
		  A_{j+.5,i+.5} & A_{j+.5,i+1} & A_{j+.5,i+1.5} & \Sigma A_{j+.5,i+.5}\\
		  A_{j+1,i+.5} & A_{j+1,i+1} & A_{j+1,i+1.5} & \Sigma A_{j+1,i+.5}\\
		  A_{j+1.5,i+.5} & A_{j+1.5,i+1} & A_{j+1.5,i+1.5} & \Sigma A_{j+1.5,i+.5}\\
		  \Sigma A_{j+.5,i+.5} & \Sigma A_{j+.5,i+1} & \Sigma A_{j+.5,i+1.5} & \Sigma^2 A_{j+.5,i+.5}\\
		};

		\path[->, font = \tiny, auto]
		  (d-1-1) edge node{$\alpha_{j+.5,i}''$} (d-1-2)
		  (d-1-2) edge[densely dotted] node{$\alpha_{j+.5,i+1}'$} (d-1-3)
		  (d-1-3) edge[densely dotted] node{$\alpha_{j+.5,n}^{i+1}$} (d-1-4)
		  
		  (d-2-1) edge node{$\alpha_{j+1,i}''$} (d-2-2)
		  (d-2-2) edge[densely dotted] node{$\alpha_{j+1,i+1}'$} (d-2-3)
		  (d-2-3) edge[densely dotted] node{$\alpha_{j+1,n}^{i+1}$} (d-2-4)
		  
		  (d-3-1) edge node{$\alpha_{j+1.5,i}''$} (d-3-2)
		  (d-3-2) edge[densely dotted] node{$\alpha_{j+1.5,i+1}'$} (d-3-3)
		  (d-3-3) edge[densely dotted] node{$\alpha_{j+1.5,n}^{i+1}$} (d-3-4)
		  
		  (d-4-1) edge node{$\Sigma\alpha_{j+.5,i}''$} (d-4-2)
		  (d-4-2) edge node{$\Sigma\alpha_{j+.5,i+1}'$} (d-4-3)
		  (d-4-3) edge node{$\Sigma\alpha_{j+.5,n}^{i+1}$} (d-4-4)
		  
		  (d-1-1) edge node{$\varphi_{j,i+.5}''$} (d-2-1)
		  (d-1-2) edge node{$\varphi_{j,i+1}''$} (d-2-2)
		  (d-1-3) edge[densely dashed] node{$\varphi_{j,i+1.5}''$} (d-2-3)
		  (d-1-4) edge node{$\Sigma\varphi_{j,i+.5}''$} (d-2-4)
		  
		  (d-2-1) edge node{$\varphi_{j+1,i+.5}'$} (d-3-1)
		  (d-2-2) edge node{$\varphi_{j+1,i+1}'$} (d-3-2)
		  (d-2-3) edge[densely dashed] node{$\varphi_{j+1,i+1.5}'$} (d-3-3)
		  (d-2-4) edge node{$\Sigma\varphi_{j+1,i+.5}'$} (d-3-4)
		  
		  (d-3-1) edge node{$\varphi_{n,i+.5}^{j+1}$} (d-4-1)
		  (d-3-2) edge node{$\varphi_{n,i+1}^{j+1}$} (d-4-2)
		  (d-3-3) edge[densely dashed] node{$\varphi_{n,i+1.5}^{j+1}$} (d-4-3)
		  (d-3-4) edge node{$\Sigma\varphi_{n,i+.5}^{j+1}$} (d-4-4)
		  ;
		  \path (d-4-3) -- (d-3-4) node[pos=.6,xshift=-.5em] (b) {$\Omega_{ji}$};
	\end{tikzpicture}
 
	\caption{The diagram $\Gamma_{ji}$.}
\end{figure}

\noindent where $\Omega_{ji}$ is anticommutative. We may call the above diagram $\Gamma_{ji}$.
For the previously fixed $i$, we now graft together the diagrams $\Gamma_{ji}$ vertically; the third row of $\Gamma_{ji}$ is the first row of $\Gamma_{j+1,i}$, and so we extend this diagram by composing the vertical morphisms in
\begin{center}
	\begin{tikzpicture}	
		\diagram{d}{3em}{3.5em}{
		  A_{j+1,i+.5} & A_{j+1,i+1} & A_{j+1,i+1.5} & \Sigma A_{j+1,i+.5}\\
		  A_{j+1.5,i+.5} & A_{j+1.5,i+1} & A_{j+1.5,i+1.5} & \Sigma A_{j+1.5,i+.5}\\
		  A_{j+2,i+.5} & A_{j+2,i+1} & A_{j+2,i+1.5} & \Sigma A_{j+2,i+.5}\\
		};

		\path[->, font = \tiny, auto]
		  (d-1-1) edge node{$\alpha_{j+1,i}''$} (d-1-2)
		  (d-1-2) edge node{$\alpha_{j+1,i+1}'$} (d-1-3)
		  (d-1-3) edge node{$\alpha_{j+1,n}^{i+1}$} (d-1-4)
		  
		  (d-2-1) edge node{$\alpha_{j+1.5,i}''$} (d-2-2)
		  (d-2-2) edge node{$\alpha_{j+1.5,i+1}'$} (d-2-3)
		  (d-2-3) edge node{$\alpha_{j+1.5,n}^{i+1}$} (d-2-4)
		  
		  (d-3-1) edge node{$\alpha_{j+2,i}''$} (d-3-2)
		  (d-3-2) edge node{$\alpha_{j+2,i+1}'$} (d-3-3)
		  (d-3-3) edge node{$\alpha_{j+2,n}^{i+1}$} (d-3-4)
		  
		  (d-1-1) edge node[pos=.4]{$\varphi_{j+1,i+.5}'$} (d-2-1)
		  (d-1-2) edge node[pos=.4]{$\varphi_{j+1,i+1}'$} (d-2-2)
		  (d-1-3) edge node[pos=.4]{$\varphi_{j+1,i+1.5}'$} (d-2-3)
		  (d-1-4) edge node[pos=.4]{$\Sigma\varphi_{j+1,i+.5}'$} (d-2-4)
		  
		  (d-2-1) edge node[pos=.4]{$\varphi_{j+1,i+.5}''$} (d-3-1)
		  (d-2-2) edge node[pos=.4]{$\varphi_{j+1,i+1}''$} (d-3-2)
		  (d-2-3) edge node[pos=.4]{$\varphi_{j+1,i+1.5}''$} (d-3-3)
		  (d-2-4) edge node[pos=.4]{$\Sigma\varphi_{j+1,i+.5}''$} (d-3-4)
		  ;
	\end{tikzpicture}
\end{center}
to get 
\[
		\varphi_{j+1,i+.5}=\varphi_{j+1,i+.5}''\varphi_{j+1,i+.5}'\quad\mt{and}\quad\varphi_{j+1,i+1}=\varphi_{j+1,i+1}''\varphi_{j+1,i+1}'
\]
vertically by construction, and 
\[
	\varphi_{j+1,i+1.5}\coloneqq\varphi_{j+1,i+1.5}''\varphi_{j+1,i+1.5}'
\]
by definition. Similarly, the suspension of the third row of $\Gamma_{ji}$ is the last row of $\Gamma_{j+1,i}$, so we complete our extension of $\Gamma_{ji}$ by composing the vertical morphisms in
\begin{center}
	\begin{tikzpicture}	
		\diagram{d}{3em}{3.5em}{
		  A_{j+2.5,i+.5} & A_{j+2.5,i+1} & A_{j+2.5,i+1.5} & \Sigma A_{j+2.5,i+.5}\\
		  \Sigma A_{j+1.5,i+.5} & \Sigma A_{j+1.5,i+1} & \Sigma A_{j+1.5,i+1.5} & \Sigma^2 A_{j+1.5,i+.5}\\
		  \Sigma^2 A_{j+.5,i+.5} & \Sigma^2 A_{j+.5,i+1} & \Sigma^2 A_{j+.5,i+1.5} & \Sigma^3 A_{j+.5,i+.5}.\\
		};

		\path[->, font = \tiny, auto]
		  (d-1-1) edge node{$\alpha_{j+2.5,i}''$} (d-1-2)
		  (d-1-2) edge node{$\alpha_{j+2.5,i+1}'$} (d-1-3)
		  (d-1-3) edge node{$\alpha_{j+2.5,n}^{i+1}$} (d-1-4)
		  
		  (d-2-1) edge node{$\Sigma\alpha_{j+1.5,i}''$} (d-2-2)
		  (d-2-2) edge node{$\Sigma\alpha_{j+1.5,i+1}'$} (d-2-3)
		  (d-2-3) edge node{$\Sigma\alpha_{j+1.5,n}^{i+1}$} (d-2-4)
		  
		  (d-3-1) edge node{$\Sigma^2\alpha_{j+.5,i}''$} (d-3-2)
		  (d-3-2) edge node{$\Sigma^2\alpha_{j+.5,i+1}'$} (d-3-3)
		  (d-3-3) edge node{$\Sigma^2Ê\alpha_{j+.5,n}^{i+1}$} (d-3-4)
		  
		  (d-1-1) edge node[pos=.4]{$\varphi_{n,i+.5}^{j+2}$} (d-2-1)
		  (d-1-2) edge node[pos=.4]{$\varphi_{n,i+1}^{j+2}$} (d-2-2)
		  (d-1-3) edge node[pos=.4]{$\varphi_{n,i+1.5}^{j+2}$} (d-2-3)
		  (d-1-4) edge node[pos=.4]{$\Sigma\varphi_{n,i+.5}^{j+2}$} (d-2-4)
		  
		  (d-2-1) edge node[pos=.4]{$\Sigma\varphi_{n,i+.5}^{j+1}$} (d-3-1)
		  (d-2-2) edge node[pos=.4]{$\Sigma\varphi_{n,i+1}^{j+1}$} (d-3-2)
		  (d-2-3) edge node[pos=.4]{$\Sigma\varphi_{n,i+1.5}^{j+1}$} (d-3-3)
		  (d-2-4) edge node[pos=.4]{$\Sigma^2\varphi_{n,i+.5}^{j+1}$} (d-3-4)
		  ;
		  \path (d-2-3) -- (d-1-4) node[pos=.6,xshift=.5em] (b) {$\Omega_{j+1,i}$};
		  \path (d-3-3) -- (d-2-4) node[pos=.6,xshift=.5em] (b) {$\Sigma\Omega_{ji}$};
	\end{tikzpicture}
\end{center}
We refrain from naming these compositions at this point, but we notice that the rectangle enveloping $\Omega_{j+1,i}$ and $\Sigma\Omega_{ji}$ is commutative, since the individual squares are anticommutative. Composing the diagram vertically, we will think of the enveloping rectangle as the vertical composition of $\Omega_{j+1,i}$ and $\Sigma\Omega_{ji}$ and thus we denote it by $\Sigma\Omega_{ji}\circ_v \Omega_{j+1,i}$. All told, we obtain a commutative diagram
\begin{center}
		\begin{tikzpicture}	
		\smalldiagram{d}{3em}{3.5em}{
		  A_{j+.5,i+.5} & A_{j+.5,i+1} & A_{j+.5,i+1.5} & \Sigma A_{j+.5,i+.5}\\
		  A_{j+1,i+.5} & A_{j+1,i+1} & A_{j+1,i+1.5} & \Sigma A_{j+1,i+.5}\\
		  A_{j+2,i+.5} & A_{j+2,i+1} & A_{j+2,i+1.5} & \Sigma A_{j+2,i+.5}\\
		  A_{j+2.5,i+.5} & A_{j+2.5,i+1} & A_{j+2.5,i+1.5} & \Sigma A_{j+2.5,i+.5}\\
		  \Sigma^2 A_{j+.5,i+.5} & \Sigma^2 A_{j+.5,i+1} & \Sigma^2 A_{j+.5,i+1.5} & \Sigma^3 A_{j+.5,i+.5}.\\
		};

		\path[->, font = \tiny, auto]
		  (d-1-1) edge node{$\alpha_{j+.5,i}''$} (d-1-2)
		  (d-1-2) edge node{$\alpha_{j+.5,i+1}'$} (d-1-3)
		  (d-1-3) edge node{$\alpha_{j+.5,n}^{i+1}$} (d-1-4)
		  
		  (d-2-1) edge node{$\alpha_{j+1,i}''$} (d-2-2)
		  (d-2-2) edge node{$\alpha_{j+1,i+1}'$} (d-2-3)
		  (d-2-3) edge node{$\alpha_{j+1,n}^{i+1}$} (d-2-4)
		  
		  (d-3-1) edge node{$\alpha_{j+2,i}''$} (d-3-2)
		  (d-3-2) edge node{$\alpha_{j+2,i+1}'$} (d-3-3)
		  (d-3-3) edge node{$\alpha_{j+2,n}^{i+1}$} (d-3-4)
		  
		  (d-4-1) edge node{$\alpha_{j+2.5,i}''$} (d-4-2)
		  (d-4-2) edge node{$\alpha_{j+2.5,i+1}'$} (d-4-3)
		  (d-4-3) edge node{$\alpha_{j+2.5,n}^{i+1}$} (d-4-4)

		  (d-5-1) edge node{$\Sigma^2\alpha_{j+.5,i}''$} (d-5-2)
		  (d-5-2) edge node{$\Sigma^2\alpha_{j+.5,i+1}'$} (d-5-3)
		  (d-5-3) edge node{$\Sigma^2\alpha_{j+.5,n}^{i+1}$} (d-5-4)
		  
		  (d-1-1) edge node{$\varphi_{j,i+.5}''$} (d-2-1)
		  (d-1-2) edge node{$\varphi_{j,i+1}''$} (d-2-2)
		  (d-1-3) edge node{$\varphi_{j,i+1.5}''$} (d-2-3)
		  (d-1-4) edge node{$\Sigma\varphi_{j,i+.5}''$} (d-2-4)
		  
		  (d-2-1) edge[densely dashed] node{$\varphi_{j+1,i+.5}$} (d-3-1)
		  (d-2-2) edge[densely dashed] node{$\varphi_{j+1,i+1}$} (d-3-2)
		  (d-2-3) edge[densely dashed] node{$\varphi_{j+1,i+1.5}$} (d-3-3)
		  (d-2-4) edge[densely dashed] node{$\Sigma\varphi_{j+1,i+.5}$} (d-3-4)
		  
		  (d-3-1) edge node{$\varphi_{j+2,i+.5}'$} (d-4-1)
		  (d-3-2) edge node{$\varphi_{j+2,i+1}'$} (d-4-2)
		  (d-3-3) edge node{$\varphi_{j+2,i+1.5}'$} (d-4-3)
		  (d-3-4) edge node{$\Sigma\varphi_{j+2,i+.5}'$} (d-4-4)

		  (d-4-1) edge[densely dashed] node{$$} (d-5-1)
		  (d-4-2) edge[densely dashed] node{$$} (d-5-2)
		  (d-4-3) edge[densely dashed] node{$$} (d-5-3)
		  (d-4-4) edge[densely dashed] node{$$} (d-5-4)
		  ;
	\end{tikzpicture}
 
\end{center}
where the dashed morphisms represent the composite maps we constructed above. Applying this specific construction to $\Gamma_{1i}$, i.e.\ the diagram (\ref{gamma1i}), we denote the resulting $5\times 4$ diagram by $\Gamma_{i}^1$.
We now repeat the procedure to $\Gamma_i^1$, noticing that the fourth row of this diagram corresponds to the first row of $\Gamma_{3i}$ and the suspension of the fourth row of $\Gamma_i^1$ is the last row of $\Gamma_{3i}$. This grafting procedure thus produces a $6\times 4$ diagram which we denote by $\Gamma_i^2$, and in which the lowest row is $\Sigma^3A_{1\bullet}$ and the lower right hand square is anticommutative since it is the one we denote by
\[
	\Sigma^2\Omega_{1i}\circ_v \Sigma\Omega_{2i}\circ_v \Omega_{3i}.
\] 
We now repeatedly graft on $\Gamma_{ji}$ to $\Gamma^{j-2}_i$ for $j\in\{3,\dots,n-2\}$ until we arrive at the $(n+1)\times 4$ diagram $\textup I_i\coloneqq\Gamma^{n-3}_i$
\begin{figure}[!h]
	\centering
	\begin{tikzpicture}	
		\diagram{d}{2em}{2em}{
		  A_{1,i+.5} & A_{1,i+1} & A_{1,i+1.5} & \Sigma A_{1,i+.5}\\
		  A_{2,i+.5} & A_{2,i+1} & A_{2,i+1.5} & \Sigma A_{2,i+.5}\\
		  A_{3,i+.5} & A_{3,i+1} & A_{3,i+1.5} & \Sigma A_{3,i+.5}\\
		  \vdots & \vdots & \vdots & \vdots\\
		  A_{n,i+.5} & A_{n,i+1} & A_{n,i+1.5} & \Sigma A_{n,i+.5}\\
		  \Sigma^{n-2} A_{1,i+.5} & \Sigma^{n-2} A_{1,i+1} & \Sigma^{n-2} A_{1,i+1.5} & \Sigma^{n-1} A_{1,i+.5}.\\
		};

		\path[->, font = \scriptsize, auto]
		  (d-1-1) edge node{$$} (d-1-2)
		  (d-1-2) edge node{$$} (d-1-3)
		  (d-1-3) edge node{$$} (d-1-4)
		  
		  (d-2-1) edge node{$$} (d-2-2)
		  (d-2-2) edge node{$$} (d-2-3)
		  (d-2-3) edge node{$$} (d-2-4)
		  
		  (d-3-1) edge node{$$} (d-3-2)
		  (d-3-2) edge node{$$} (d-3-3)
		  (d-3-3) edge node{$$} (d-3-4)
		  
		  (d-5-1) edge node{$$} (d-5-2)
		  (d-5-2) edge node{$$} (d-5-3)
		  (d-5-3) edge node{$$} (d-5-4)

		  (d-6-1) edge node{$$} (d-6-2)
		  (d-6-2) edge node{$$} (d-6-3)
		  (d-6-3) edge node{$$} (d-6-4)
		  
		  (d-1-1) edge node{$$} (d-2-1)
		  (d-1-2) edge node{$$} (d-2-2)
		  (d-1-3) edge node{$$} (d-2-3)
		  (d-1-4) edge node{$$} (d-2-4)
		  
		  (d-2-1) edge node{$$} (d-3-1)
		  (d-2-2) edge node{$$} (d-3-2)
		  (d-2-3) edge node{$$} (d-3-3)
		  (d-2-4) edge node{$$} (d-3-4)
		  
		  (d-3-1) edge node{$$} (d-4-1)
		  (d-3-2) edge node{$$} (d-4-2)
		  (d-3-3) edge node{$$} (d-4-3)
		  (d-3-4) edge node{$$} (d-4-4)

		  (d-4-1) edge node{$$} (d-5-1)
		  (d-4-2) edge node{$$} (d-5-2)
		  (d-4-3) edge node{$$} (d-5-3)
		  (d-4-4) edge node{$$} (d-5-4)

		  (d-5-1) edge node{$$} (d-6-1)
		  (d-5-2) edge node{$$} (d-6-2)
		  (d-5-3) edge node{$$} (d-6-3)
		  (d-5-4) edge node{$$} (d-6-4)
		  ;
		  \path (d-6-3) -- (d-5-4) node[midway] (b) {$\Omega_{i}$};
	\end{tikzpicture}
	\caption{The diagram $\textup I_i$.}
\end{figure}
where $\Omega_i$ represents the square
\[
	\Sigma^{n-3}\Omega_{1i}\circ_v \Sigma^{n-4}\Omega_{2i}\circ_v \cdots\circ_v\Sigma\Omega_{n-3,i}\circ_v \Omega_{n-2,i}.
\]
Note that $\Omega_i$ is a composition of $n-2$ anticommutative squares, so $\Omega_i$ commutes if $n$ is even and anticommutes if $n$ is odd. Note also that the first three columns are $n$-angles in the cluster tilting subcategory, while the first $n$ rows of $\textup I_i$ are distinguished triangles in the ambient triangulated category. \smallskip

We now essentially repeat this procedure in the horizontal direction, grafting $\textup I_i$ and $\textup I_{i+1}$ for $i\in\{1,\dots,n-3\}$. This means that we notice that the third column of $\textup I_1$ is the first column of $\textup I_2$ and the suspension of the third column of $\textup I_1$ is the last column of $\textup I_2$, and so a similar grafting procedure takes place to create the $(n+1)\times 5$ diagram
\begin{center}
	\begin{tikzpicture}	
		\diagram{d}{2em}{1.2em}{
		  A_{11} & A_{12} & A_{13} & A_{1,3.5} & \Sigma^2 A_{11}\\
		  A_{21} & A_{22} & A_{23} & A_{2,3.5} & \Sigma^2 A_{21}\\
		  A_{31} & A_{32} & A_{33} & A_{3,3.5} & \Sigma^2 A_{31}\\
		  \vdots & \vdots & \vdots & \vdots &\vdots\\
		  A_{n1} & A_{n2} & A_{n3} & A_{n,3.5} & \Sigma^2 A_{n1}\\
		  \Sigma^{n-2} A_{11} & \Sigma^{n-2} A_{12} & \Sigma^{n-2} A_{13} & \Sigma^{n-2}A_{1,3.5} & \Sigma^{n} A_{11}.\\
		};

		\path[->, font = \scriptsize, auto]
		  (d-1-1) edge node{$$} (d-1-2)
		  (d-1-2) edge node{$$} (d-1-3)
		  (d-1-3) edge node{$$} (d-1-4)
		  (d-1-4) edge node{$$} (d-1-5)
		  
		  (d-2-1) edge node{$$} (d-2-2)
		  (d-2-2) edge node{$$} (d-2-3)
		  (d-2-3) edge node{$$} (d-2-4)
		  (d-2-4) edge node{$$} (d-2-5)
		  
		  (d-3-1) edge node{$$} (d-3-2)
		  (d-3-2) edge node{$$} (d-3-3)
		  (d-3-3) edge node{$$} (d-3-4)
		  (d-3-4) edge node{$$} (d-3-5)
		  
		  (d-5-1) edge node{$$} (d-5-2)
		  (d-5-2) edge node{$$} (d-5-3)
		  (d-5-3) edge node{$$} (d-5-4)
		  (d-5-4) edge node{$$} (d-5-5)
		  
		  (d-6-1) edge node{$$} (d-6-2)
		  (d-6-2) edge node{$$} (d-6-3)
		  (d-6-3) edge node{$$} (d-6-4)
		  (d-6-4) edge node{$$} (d-6-5)
		  
		  (d-1-1) edge node{$$} (d-2-1)
		  (d-1-2) edge node{$$} (d-2-2)
		  (d-1-3) edge node{$$} (d-2-3)
		  (d-1-4) edge node{$$} (d-2-4)
		  (d-1-5) edge node{$$} (d-2-5)
		  
		  (d-2-1) edge node{$$} (d-3-1)
		  (d-2-2) edge node{$$} (d-3-2)
		  (d-2-3) edge node{$$} (d-3-3)
		  (d-2-4) edge node{$$} (d-3-4)
		  (d-2-5) edge node{$$} (d-3-5)
		  
		  (d-3-1) edge node{$$} (d-4-1)
		  (d-3-2) edge node{$$} (d-4-2)
		  (d-3-3) edge node{$$} (d-4-3)
		  (d-3-4) edge node{$$} (d-4-4)
		  (d-3-5) edge node{$$} (d-4-5)
		  
		  (d-4-1) edge node{$$} (d-5-1)
		  (d-4-2) edge node{$$} (d-5-2)
		  (d-4-3) edge node{$$} (d-5-3)
		  (d-4-4) edge node{$$} (d-5-4)
		  (d-4-5) edge node{$$} (d-5-5)
		  
		  (d-5-1) edge node{$$} (d-6-1)
		  (d-5-2) edge node{$$} (d-6-2)
		  (d-5-3) edge node{$$} (d-6-3)
		  (d-5-4) edge node{$$} (d-6-4)
		  (d-5-5) edge node{$$} (d-6-5)
		  ;
		  \path (d-6-4) -- (d-5-5) node[midway] (b) {$\Sigma\Omega_1\circ_h\Omega_2$};
	\end{tikzpicture}
\end{center}
 which we will denote by $\textup I_1^1$. Note that $\Omega_1$ and $\Omega_2$ either both commute or both anticommute, so the enveloping square, which we denote as the horisontal composition $\Sigma\Omega_1\circ_h\Omega_2$, commutes regardless of the parity of $n$. Continuing on in the same fashion to before, we end with the $(n+1)\times(n+1)$ diagram 
\begin{center}
	\begin{tikzpicture}	
		\diagram{d}{2em}{1.2em}{
		  A_{11} & A_{12} & A_{13} & \cdots & A_{1n} & \Sigma^{n-2} A_{11}\\
		  A_{21} & A_{22} & A_{23} & \cdots & A_{2n} & \Sigma^{n-2} A_{21}\\
		  A_{31} & A_{32} & A_{33} & \cdots & A_{3n} & \Sigma^{n-2} A_{31}\\
		  \vdots & \vdots & \vdots & &  \vdots &\vdots\\
		  A_{n1} & A_{n2} & A_{n3} & \cdots & A_{n,n} & \Sigma^{n-2} A_{n1}\\
		  \Sigma^{n-2} A_{11} & \Sigma^{n-2} A_{12} & \Sigma^{n-2} A_{13} & \cdots & \Sigma^{n-2}A_{1n} & \Sigma^{2(n-2)} A_{11},\\
		};

		\path[->, font = \scriptsize, auto]
		  (d-1-1) edge node{$$} (d-1-2)
		  (d-1-2) edge node{$$} (d-1-3)
		  (d-1-3) edge node{$$} (d-1-4)
		  (d-1-4) edge node{$$} (d-1-5)
		  (d-1-5) edge node{$$} (d-1-6)
		  
		  (d-2-1) edge node{$$} (d-2-2)
		  (d-2-2) edge node{$$} (d-2-3)
		  (d-2-3) edge node{$$} (d-2-4)
		  (d-2-4) edge node{$$} (d-2-5)
		  (d-2-5) edge node{$$} (d-2-6)
		  
		  (d-3-1) edge node{$$} (d-3-2)
		  (d-3-2) edge node{$$} (d-3-3)
		  (d-3-3) edge node{$$} (d-3-4)
		  (d-3-4) edge node{$$} (d-3-5)
		  (d-3-5) edge node{$$} (d-3-6)
		  
		  (d-5-1) edge node{$$} (d-5-2)
		  (d-5-2) edge node{$$} (d-5-3)
		  (d-5-3) edge node{$$} (d-5-4)
		  (d-5-4) edge node{$$} (d-5-5)
		  (d-5-5) edge node{$$} (d-5-6)
		  
		  (d-6-1) edge node{$$} (d-6-2)
		  (d-6-2) edge node{$$} (d-6-3)
		  (d-6-3) edge node{$$} (d-6-4)
		  (d-6-4) edge node{$$} (d-6-5)
		  (d-6-5) edge node{$$} (d-6-6)
		  
		  (d-1-1) edge node{$$} (d-2-1)
		  (d-1-2) edge node{$$} (d-2-2)
		  (d-1-3) edge node{$$} (d-2-3)
		  (d-1-5) edge node{$$} (d-2-5)
		  (d-1-6) edge node{$$} (d-2-6)
		  
		  (d-2-1) edge node{$$} (d-3-1)
		  (d-2-2) edge node{$$} (d-3-2)
		  (d-2-3) edge node{$$} (d-3-3)
		  (d-2-5) edge node{$$} (d-3-5)
		  (d-2-6) edge node{$$} (d-3-6)
		  
		  (d-3-1) edge node{$$} (d-4-1)
		  (d-3-2) edge node{$$} (d-4-2)
		  (d-3-3) edge node{$$} (d-4-3)
		  (d-3-5) edge node{$$} (d-4-5)
		  (d-3-6) edge node{$$} (d-4-6)
		  
		  (d-4-1) edge node{$$} (d-5-1)
		  (d-4-2) edge node{$$} (d-5-2)
		  (d-4-3) edge node{$$} (d-5-3)
		  (d-4-5) edge node{$$} (d-5-5)
		  (d-4-6) edge node{$$} (d-5-6)
		  
		  (d-5-1) edge node{$$} (d-6-1)
		  (d-5-2) edge node{$$} (d-6-2)
		  (d-5-3) edge node{$$} (d-6-3)
		  (d-5-5) edge node{$$} (d-6-5)
		  (d-5-6) edge node{$$} (d-6-6)
		  ;
		  \path (d-6-5) -- (d-5-6) node[midway] (b) {$\Omega$};
	\end{tikzpicture}
\end{center}
denoted $\textup I\coloneqq\textup I_1^{n-3}$, such that the first $n$ rows and columns are all $n$-angles and such that the upper most $2\times(n+1)$ diagram is the morphism we started with. Note, lastly, that the square $\Omega$ represents the square
\[
	\Sigma^{n-3}\Omega_1\circ_h\Sigma^{n-4}\Omega_2\circ_h\cdots\circ_h\Sigma\Omega_{n-3}\circ_h\Omega_{n-2}.
\]
which is a composition of $n-2$ squares that either all commute or all anticommute. As such, $\Omega$ commutes exactly when $n$ is even and anticommutes when $n$ is odd.
\end{prf}


\section{Further explorations}

For this section, let $R$ be a commutative local ring with principal non-zero maximal ideal $\mathfrak m=(p)$ satisfying $\frak m^2=0$. It was originally proved in \cite{bjt} that the category, \cat C, of finitely generated free $R$-modules can be equipped with an $n$-angulated structure for $n\>3$ if $n$ is even, and also for $n$ odd if $2p=0$ such that the suspension is the identity functor on \cat C. An example of an $n$-angle in \cat C is 
\begin{center}
		\begin{tikzpicture}	
			\diagram{d}{2em}{3em}{
			  R & R & R & \cdots & R & \Sigma R\\ 
			};
	
			\path[->, font = \scriptsize, auto]
			  (d-1-1) edge node{$p$} (d-1-2)
			  (d-1-2) edge node{$p$} (d-1-3)
			  (d-1-3) edge node{$p$} (d-1-4)
			  (d-1-4) edge node{$p$} (d-1-5)
			  (d-1-5) edge node{$p$} (d-1-6)
			  ;
		\end{tikzpicture}
	\end{center}
which is commonly denoted $R(p)_\bullet$. We suspect that for arbitrary $n\>3$, the morphism of $n$-angles
\[
	(0,p,0,\dots,0)\colon R(p)_\bullet \longrightarrow R(p)_\bullet
\]
does not extend to a middling good diagram. However, the author has been unable to prove or disprove this conjecture and, so, we leave it as an open question. It also remains to be seen whether all good and Verdier good morphisms of $n$-angles are middling good, as is the case in triangulated categories.

\paragraph{Questions:} 
\begin{enumerate}
\item[$\bullet$] Does there exist an $n$-angulated category and a morphism of $n$-angles in said category that does not extend to a middling good diagram?
\item[$\bullet$] Is it always the case that good morphisms are middling good?
\item[$\bullet$] Is it always the case that Verdier good morphisms are middling good?
\end{enumerate}

\printbibliography

\vspace{.5cm}
\noindent{\scshape Department of Mathematical Sciences, NTNU\\ N-7491 Trondheim, Norway}\medskip

\noindent \texttt{sebastian.martensen@ntnu.no}\smallskip

\noindent{ORCID iD:} \texttt{0000-0002-4477-8856}

\end{document}